\documentclass[12pt,reqno]{amsart}
\usepackage{amssymb}
\usepackage{epsf,epsfig,amsfonts}
\usepackage{mdwlist}
\usepackage{amsfonts,epsfig,a4wide}
\usepackage[usenames,dvips]{color}

\usepackage[all]{xy}

\usepackage[latin1]{inputenc}

\usepackage{graphicx}

\numberwithin{equation}{section}

\DeclareFontFamily{U}{rsf}{\skewchar\font'177}%
\DeclareFontShape{U}{rsf}{m}{n}{<-6>rsfs5<6-8>rsfs7<8->rsfs10}{}%
\DeclareFontShape{U}{rsf}{b}{n}{<-6>rsfs5<6-8>rsfs7<8->rsfs10}{}%
\DeclareMathAlphabet\RSFS{U}{rsf}{m}{n}
\SetMathAlphabet\RSFS{bold}{U}{rsf}{b}{n}
\let\scr=\rfs

\newtheorem{thm}{Theorem}[section]
\newtheorem{lem}{Lemma}[section]
\newtheorem{prop}{Proposition}[section]
\newtheorem{defi}{Definition}[section]
\newtheorem{corol}{Corollary}[section]
\newtheorem{rem}{Remark}[section]
\newcommand{\gG}{g}

\def\sf#1{{\mathsf{#1}}} 

\def\slsf#1{{\slshape \sffamily #1\/}}
\def\xrar#1{\xrightarrow{#1}}

\def\pmatr#1{{\left(\begin{matrix} #1 \end{matrix} \right)}}
\def\smpmatr#1{{\left(\begin{smallmatrix} #1 \end{smallmatrix}\right)}}

\def\state#1. {\noindent{\bf{#1.} }}

\def\cc{{\mathbb C}} 
\def\rr{{\mathbb R}}

\def\bfk{{\mathbf{K}}}

\def\bft{{\mathbf{T}}}

\def\cala{{\mathcal{A}}}

\def\cald{{\mathcal{D}}}
\def\cale{{\mathcal{E}}}
  
\def\calf{{\mathcal{F}}}
\def\calg{{\mathcal{G}}}

\def\calk{{\mathcal{K}}}

\def\scrl{{\scr{L}}}
\def\cals{{\mathcal{S}}}

\def\d{\partial}
\def\half{{\frac12}}

\def\vph{^{\mathstrut}}

\def\barz{{\bar z}}
\let\bz=\barz
\def\barp{{\bar p}}
\let\wt=\widetilde

\def\const{{\sf{const}}}

\def\cplxi{{\mskip2mu\sf{i}\mskip2mu}}

\def\id{{\sf{Id}}}
\def\Div{{\sf{div}}}  
\def\det{{\sf{det}}}  
\def\sym{{\sf{Sym}}}

\def\min{{\sf{min}}}

\def\log{{\sf{log}}}
\def\lim{{\sf{lim}}}

\def\tr{{\sf{Tr}}}

\def\inv{^{-1}}

\def\leq{\leqslant}
\def\geq{\geqslant}

\def\sli{{\sl i)\ \;}}
       \def\sliip{{\sl i$_{\!}$i)}}\def\slii{{\sliip\ \;}}
\def\sliiip{{\sl i$_{\!}$i$_{\!}$i)}}
       \def\sliii{{\sliiip\ \;}}

\newcommand{\weg}[1]{}

{\catcode`\@=11
\gdef\n@te#1#2{\leavevmode\vadjust{%
 {\setbox\z@\hbox to\z@{\strut#1}%
  \setbox\z@\hbox{\raise\dp\strutbox\box\z@}\ht\z@=\z@\dp\z@=\z@%
  #2\box\z@}}}
\gdef\leftnote#1{\n@te{\hss#1\quad}{}}
\gdef\rightnote#1{\n@te{\quad\kern-\leftskip#1\hss}{\moveright\hsize}}
\gdef\?{\FN@\qumark}
\gdef\qumark{\ifx\next"\DN@"##1"{\leftnote{\rm##1}}\else
 \DN@{\leftnote{\rm??}}\fi{\rm??}\next@}}

\definecolor{DarkGreen}{rgb}{0,0.6,0}

\definecolor{MyMagenta}{rgb}{0.6,0.0,0.5}

\begin{document}

\baselineskip =14.7pt plus 2pt
\leftmargini=27pt 
\leftmarginii=21pt 
\leftmarginiii=19pt 
\leftmarginiv=19pt  

\title[Cubic  integrals for geodesic flows]%
{Differential invariants for cubic integrals  of  geodesic flows on surfaces }
\author[V.~S.~Matveev]{Vladimir  S.~Matveev}
\address{Mathematisches Institut\\ 
 Fakultät für Mathematik und Informatik\\ 
 Friedrich-Schiller-Universität Jena\\ 
 07737  Jena } 
\email{vladimir.matveev@uni-jena.de}
\author[V.~V.~Shevchishin]{Vsevolod V.~Shevchishin}
\address{Mathematisches Institut\\ 
 Fakultät für Mathematik und Informatik\\ 
 Friedrich-Schiller-Universität Jena\\ 
 07737  Jena 
} \email{shevchishin@googlemail.com }


\begin{abstract}
 We construct  differential invariants that vanish if and only if the
 geodesic flow of a 2-dimensional metric admits an integral of 3rd degree in
 momenta with a given Birkhoff-Kolokoltsov 3-codifferential.
\end{abstract}
\maketitle

\section{Introduction}
\subsection{Definitions and results} \label{intr1} Let $S$ be a surface (i.e.,
2-dimensional real manifold) equipped with a Riemannian metric $g$ given in
local coordinates by $g=\sum_{ij}g_{ij}dx^idx^j$.  Since the metric $g$ allows
us to identify the tangent and cotangent bundles of $S$, we have a scalar
product and a norm on every cotangent plane. The \slsf{geodesic flow} of the
metric $\gG$ is the Hamiltonian system with the Hamiltonian
$H:=\frac{1}{2}|\vec p|^2= \frac{1}{2} g^{ij}p_ip_j$, where
$\vec p=(p_1,p_2)$ are the momenta and $|.|$ is the norm induced by $g$.

We say that  a function $F:T^*S\to \mathbb{R}$ is an \slsf{integral}  of the
 geodesic flow of $g$  \slsf{cubic in momenta}
  (shortly: \slsf{cubic integral} for the  metric $g$), if 
\begin{enumerate} 
\item in the local coordinates $x:=x^1,y:=x^2$ on the surface and the
 corresponding momenta $p_x,p_y$, $F$ is a homogeneous polynomial in the
 momenta of degree 3: 
\begin{equation} \label{eqn:int}
  F(x,y;p_x,p_y)= a_0(x,y) p_x^3 + a_1(x,y) p_x^2p_y + 
a_2(x,y) p_x p_y^2 + a_3(x,y)p_y^3, 
\end{equation}

\item $F$ is an integral of the geodesic flow of $g$, i.e., $\{H,F\}=0$, where
 $\{\cdot, \cdot\}$ is the canonical Poisson bracket on $T^*S$.
\end{enumerate} 

Cubic integrals allow to construct different geometric structures on the
surface; we shall use the so called \slsf{holomorphic Birkhoff-Kolokoltsov
 3-codifferential}.  We consider the complex structure on $S$ given by $g$: a local
complex coordinate $z=x + \cplxi y$ is determined by the property that the
metric $g$ has the isothermic form 
\begin{equation} \label{eqn:iso} 
g =  \lambda(x,y) (dx^2 + dy^2) = \lambda(z, \bar z) dz  d\bar z.
\end{equation}
Then, the Birkhoff-Kolokoltsov 3-codifferential (corresponding to the cubic
integral \eqref{eqn:int}) is a complex-valued symmetric $(3,0)$-tensor given by
\begin{equation} \label{eqn:kol}
\textstyle
 A(z):= \Big(\left(a_0(x,y)-a_2(x,y)\right) + \cplxi\left(a_1(x,y)-a_3(x,y)\right)\Big) 
\frac{\partial}{\partial z} \otimes  \frac{\partial}{\partial z} \otimes \frac{\partial}{\partial z}.
\end{equation}
As it was shown by Birkhoff \cite{Bi}, the coefficient %
$(a_0-a_2) + \cplxi(a_1-a_3)$ is a holomorphic function of the complex
coordinate $z=x + \cplxi y$. Kolokoltsov observed in \cite{Kol} that 
under a holomorphic (\,$\Leftrightarrow$ conformal) coordinate
change this coefficient changes as the corresponding coefficient of a complex
$(3,0)$-tensor. Below in §\,\ref{cplx-calc} 
we shall prove both properties of $A$ and give several
useful formulas for $A$ in arbitrary (not necessarily isothermic)
coordinates, such as how to compute   $A$ from $F$  and how to check the 
holomorphicity of $A$, see Lemma \ref{A-form}.

The goal of this article is to give an
answer to the following question:

\begin{itemize}
\item[\ ] {\it Given a metric $g$ and a holomorphic 3-codifferential $A$, how
  to decide whether there\linebreak exists a cubic integral for the geodesic
  flow of $g$ whose Birkhoff-Kolokoltsov 3-codif­ferential coincides with $A$?}
\end{itemize}

If such a cubic integral $F$ exists, we say that \slsf{$A$ is compatible with
 $g$}. 

A \emph{complete algorithmic answer} is given in §\,\ref{algor} using Theorems
\ref{th:main1}, \ref{th:main2}{ and also Proposi­tions
 \ref{R=const}, \ref{A-kill}}. Moreover,  in the most interesting case,
 Theorems \ref{th:main1}, \ref{th:main2} give an explicit formula for the
 cubic integral.

\subsection{Main theorems.} 
In the formulas below we use Einstein's summation convention, ie., we sum over
repeating upper and lower indices. The notation $a_{;jk}$ and so on means
covariant derivation(s) of a function (or a tensor) $a$,  and 
$a^{(ijk)}$ denotes the symmetrisation in indices $i,j,k$.

Let $x=x^1,y=x^2$ be arbitrary coordinates on the surface $S$. Given a metric
$g=g_{ij}dx^jdx^j$ on $S$, let $R^i_{ijk} $ be its Riemannian curvature, and
set $\lambda(x,y):=\sqrt{\det(g_{ij})}$.  Since our considerations are local, we
equip $S$ with the orientation given by the metric volume  
form $\omega_g=\lambda dx\land dy$. 

We denote by $\{f,h\}_g:=\tfrac{1}{\lambda} \left(\tfrac{\partial f}{\partial x}\tfrac{\partial h}{\partial y}-
 \tfrac{\partial h}{\partial x}\tfrac{\partial f}{\partial y}\right)$ the Poisson bracket of functions
$f,h$ on $S$ with respect to $\omega_g$.  Further, the tensor
$J^i{}_j:=g^{ik}\omega_{kj}$ is the operator of rotation by $90^\circ$, it defines the
complex structure on $S$ compatible with the metric $g$ and the orientation,
i.e., corresponds to the multiplication by $\cplxi$ in the complex coordinate
$z=x+ \cplxi y$.

\smallskip

Define the following smooth functions $\varphi_0,\ldots,\varphi_3$, which are invariant
 algebraic expressions of the components of $g$ and its derivatives:

\begin{itemize}
\item The Gauss curvature, which  is the half of scalar curvature: 
\begin{equation}\label{phi0}
\textstyle
\varphi_0:= R:=\frac12  R^i_{jik} g^{jk};
\end{equation}
\item Half-square of the gradient of $\varphi_0=R$, half-square of the gradient of
 the result:
\begin{equation}\label{phi1phi3}\textstyle
\varphi_1:=\frac12|\nabla\varphi_0|^2 =
\frac12 g^{ij} \frac{\partial\varphi_0}{\partial x^i}\frac{\partial\varphi_0 }{\partial x^j};
\qquad 
\varphi_3:=\frac12|\nabla\varphi_1|^2 =
\frac12 g^{ij} \frac{\partial\varphi_1}{\partial x^i}\frac{\partial\varphi_1 }{\partial
  x^j}
\end{equation}
\item The Poisson bracket of $\varphi_0,\varphi_1$ with respect to $\omega_g$:
\begin{equation}\label{phi2}\textstyle
\varphi_2:= \{ \varphi_0,\varphi_1\}_g:=\frac{1}{\lambda} \left(  
 \frac{\partial\varphi_0 }{\partial x } \frac{\partial\varphi_1 }{\partial y } - \frac{\partial\varphi_0 }{\partial y } \frac{\partial
  \varphi_1 }{\partial x }     \right). 
\end{equation}
\end{itemize}

Further, let $A$ be a complex 3-codifferential on $S$ given in \emph{some}
(not necessarily isothermic) coordinates $(x^1,x^2)$ as symmetric
$(3,0)$-tensor $A^{ijk}$ with complex components.  Define $\hat A:=\Re(A)$ to be
its real part, or $\hat{A}^{ijk}:=\frac{1}{2}({A}^{ijk}+\bar{A}^{ijk})$ on the
level of components where $\bar a$ means complex conjugation. Then
$\hat{A}=(\hat{A}^{ijk})$ is a (usual) symmetric $(3,0)$-tensor with real
components.  The formula for the imaginary part $\Im(A^{ijk})$ is given in
Section \ref{invar}.

In addition to $\varphi_0,\ldots,\varphi_3$, we define the functions $D_0,\ldots,D_3$,
$\calg_2,\calg_3$, $1$-form $\calk=\calk_idx^i$, and symmetric $(3,0)$-tensors
$\hat B^{ijk},F^{ijk}$. They all are certain invariant algebraic  expressions 
in $g$,  $\hat A$, and their derivatives.
{\begin{align}
D_0:=& \textstyle
4\,\Re(A^{ijk}{}_{;ijk})= 
4\hat A^{ijk}{}_{;ijk};
\label{D0}\\
D_1:=&  g^{ij}(D_0)_{;i} (\varphi_0)_{;j} 
  -4\, \hat A^{ijk}{}_{;i}\cdot(\varphi_0)_{;j} \cdot(\varphi_0)_{;k};
\label{D1}\\
D_2:=& \{ D_0,\varphi_1\}_g + \{ \varphi_0,D_1\}_g;
\label{D2}\\
D_3:=&   g^{ij}(D_1)_{;i} (\varphi_1)_{;j} 
-4\,\hat  A^{ijk}{}_{;i}\cdot(\varphi_1)_{;j} \cdot(\varphi_1)_{;k}
\label{D3}
\end{align}} 
All these functions have a clear geometric sense: \\[-18pt]
\begin{itemize}
\item $D_0$ is the triple divergence of $4\cdot\hat A$;
\item $g^{ij}(D_0)_{;i} (\varphi_0)_{;j}$ is the scalar product $\langle\nabla D_0,\nabla\varphi_0\rangle$ of
 gradients, 
\item $\hat A^{ijk}{}_{;i}(\varphi_0)_{;j}(\varphi_0)_{;k}$ is the evaluation of
 the quadratic form corresponding to symmetric $(2,0)$-tensor
 $(\Div\hat A)^{jk}:= \hat A^{ijk}{}_{;i}$ on the $1$-form
 $d\varphi_0=dR$;
\item We obtain $\varphi_3$ from $\varphi_1$ by the same formula as $\varphi_1$ from $\varphi_0$,
 similarly, we obtain $D_3$ from $D_1$ and $\varphi_1$ by the same formula as $D_1$
 from $D_0$ and $\varphi_0$.
\end{itemize}  

Finally, define 
{\small\begin{equation}\label{G2G3}
\calg_2 :=\frac{1}{\lambda}\det
\pmatr{(\varphi_0)_{;x} & (\varphi_0)_{;y} & D_0\\
 (\varphi_1)_{;x} & (\varphi_1)_{;y} & D_1\\
 (\varphi_3)_{;x} & (\varphi_3)_{;y} & D_3 } 
\qquad  
\calg_3 :=  \frac{1}{\lambda}\, \det
\pmatr{(\varphi_0)_{;x} & (\varphi_0)_{;y} & D_0\\
 (\varphi_1)_{;x} & (\varphi_1)_{;y} & D_1\\
 (\varphi_2)_{;x} & (\varphi_2)_{;y} & D_2 }, 
\end{equation}}
\vspace{3pt}
{\small\begin{equation}\label{KBF}
\calk_i :=\frac{1}{\varphi_2}\det
\pmatr{(\varphi_0)_{;i} & D_0\\
       (\varphi_1)_{;i} & D_1} 
\qquad
\hat B^{ijk}:= g^{(ij}g^{k)l}J^{m}{}_l\calk_m
\qquad
F^{ijk} := \hat A^{ijl} +\hat B^{ijk}.
\end{equation}}

\begin{thm} \label{th:main1}%
 In the above notation, assume that the $3$-codifferential $A$ is holomorphic
 and that $\varphi_2$ in non-vanishing on $S$.  Then there exists a cubic integral
 $F$ with the Birkhoff-Kolokoltsov $3$-codifferential $A$ if and only if $g$
 and $A$ satisfy the equations $\calg_2=0,\;\calg_3=0$.
 
 Moreover, such $F$ is unique and in local coordinates $(x^1,x^2)$ 
 with the dual momenta $(p_1,p_2)$ the integral $F$ is given by
\[
F(x^1,x^2;p_1,p_2):=F^{ijk}(x^1,x^2)p_ip_jp_j
\]
where the symmetric $(3,0)$-tensor $F^{ijk}$ is computed using the formulas
\eqref{phi0}--\eqref{KBF}.
\end{thm}

\state Remarks.~1.  Our definitions and formulas involve the orientation on
the surface $S$, however, the final result is, of course, independent of the
orientation. More precisely, after reversing the orientation every real
expression (formula) either remains unchanged, or inverts the sign.
\\[3pt]
{\bf2.}\, After appropriate ``cosmetic'' changes, Theorems \ref{th:main1},
\ref{th:main2} remain valid also in the pseudo-Riemannian case, see
Proposition \ref{algo-psR}.

Let us now consider the degenerate case when $\varphi_2=\{R, \half|\nabla R|^2\}$ vanishes
identically.  In this case we set $\varphi^*_0:=\varphi_0, D^*_0:=D_0$
\begin{equation}\label{phi1*-intro}
\varphi^*_1:=\Delta R \qquad\qquad
D^*_1:= \Delta D_0  - 2\Re\big((A_{;z}R_{;z})_{;z}\big).
\end{equation}
and define the $*$-versions $\varphi^*_2,\varphi^*_3,\calg^*_2,D^*_2,D^*_3$,
$\calg^*_3,\calk^*_1,\calk^*_2...$ of the previous expressions by the
 replacing $\varphi_0, D_0,\varphi_1,...$ by $\varphi_0^*$, $D_0^*,\varphi_1^*,...$ in the formulas
 \eqref{phi2}, \eqref{D2}, \eqref{D3}, \eqref{G2G3}, \eqref{KBF} above, see
 also Section \ref{degen}. Further, let $\cald$
 be that of two expressions $\cald_i:=\det \smpmatr{(\varphi_0)_{;x^i} \;& D_0\\
 (\varphi_1)_{;x^i}\;&D_1}$ (non-$*$-versions!) %
 for which $(\varphi_0)_{;x^i}$ is non-vanishing (any of two if both
 $(\varphi_0)_{;x},(\varphi_0)_{;y}$ are non-zero.) 

\smallskip

\begin{thm} \label{th:main2} Let $g$ be a metric and $A$ a holomorphic
 $3$-codifferential.  Assume that the function $\varphi_2= \{R,\half|\nabla R|_g^2\}_g$
 vanishes identically  and the function
 $\varphi_2^*=\{R,\Delta_gR\}_g$ is non-vanishing.

 Then $g$ admits a cubic integral $F$ with the Birkhoff-Kolokoltsov
 $3$-codifferential $A$ if and only if $g$ satisfies the PDEs
 $\calg^*_2,\calg^*_3$, and $\cald$. Moreover, such $F$ is unique and in local
 coordinates $(x^1,x^2)$ with the dual momenta $(p_1,p_2)$ the integral $F$ is
 given by
$$
F(x^1,x^2;p_1,p_2):=F^{ijk}(x^1,x^2)p_ip_jp_j
$$
where the symmetric $(3,0)$-tensor $F^{ijk}$ is computed by  the $\ast-$version  
of the formula \eqref{KBF}, namely,
{\small\begin{equation}\label{KBF*}
\calk^*_i :=\frac{1}{\varphi^*_2}\det
\pmatr{(\varphi_0)_{;i} & D_0\\
       (\varphi^*_1)_{;i} & D^*_1} 
\qquad
\hat B^{ijk}:= g^{(ij}g^{k)l}J^{m}{}_l\calk^*_m
\qquad
F^{ijk} := \hat A^{ijl} +\hat B^{ijk}.
\end{equation}}
\end{thm} 

\state Remark. The equations $\calg^*_2,\calg^*_3$ are (always) necessary
conditions for compatibility of $A$ with $g$, but \emph{not sufficient} in
general. Nevertheless, if $A$ is compatible with $g$ and $\varphi^*_2$ is
non-vanishing, then the formulas \eqref{KBF*} are valid and give the correct
value of the cubic integral $F$, even in the case when $\varphi_2$ is also non-zero.

The only two cases which are not covered by Theorems \ref{th:main1},
\ref{th:main2} are the case $\varphi_0\equiv \const$ (i.e., the metric has constant
curvature) and the case $\varphi_2\equiv 0$, $\varphi_2^*\equiv 0$. These cases are easier, and will
be solved in Propositions \ref{R=const}, \ref{A-kill}.

\subsection{Algorithm for checking the existence of  cubic integrals.} \label{algor} Let $g$ be
a metric on a surface $S$ and $A$ a holomorphic $3$-codifferential.  Theorems
\ref{th:main1}, \ref{th:main2}, classical results of Darboux and Eisenhart
which we recall in \S \ref{previ}, and also Propositions \ref{R=const},
\ref{A-kill} give us the following algorithm to decide whether there exists a
cubic integral $F$ with the given Birkhoff-Kolokoltsov $3$-codifferential
$A$. Moreover, in the most interesting cases covered by Theorems
\ref{th:main1}, \ref{th:main2} we obtain an explicit formula for $F$ as
algebraic expression in $g,A$ and their derivatives. In other cases covered by
Propositions \ref{R=const}, \ref{A-kill}, the formula for $F$ requires a
solution of linear systems of ODEs.  Corresponding calculations can be easily
realized using computing software such as Maple\textsuperscript{®} or
Mathematica\textsuperscript{®}.

\vspace{0pt plus 20pt}
\pagebreak[2]
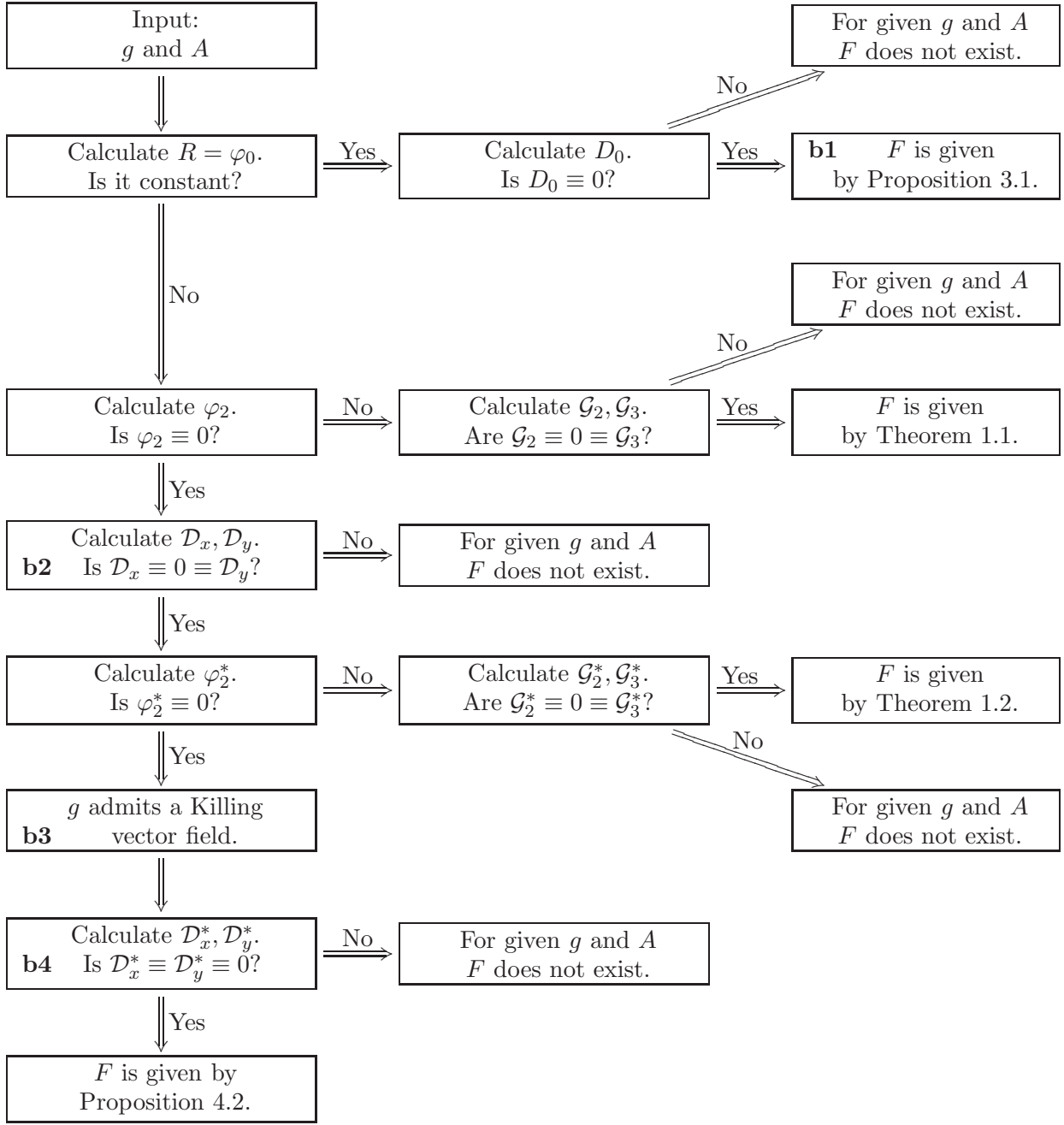
\begin{minipage}{\textwidth}
\def\algbox[#1]#2{\fbox{
\begin{minipage}{#1\textwidth}
#2 
\end{minipage}}}
\def\arrbox[#1]#2{
\begin{minipage}{#1\textwidth}
#2 
\end{minipage}}
\def\lrat{.27} 
\def\rrat{.23} 

{\small\[
\xymatrix@C+6pt{
\algbox[\lrat]{\center Input:\\ $g$ and $A$}
  \ar@{=>}[d] 
& & 
\algbox[\rrat]{\center For given $g$ and $A$\\ $F$ does not exist.}
\\
\algbox[\lrat]{\center Calculate $R=\varphi_0$.\\
Is it constant?}
\ar@{=>}[r]^{\txt{Yes}}   \ar@{=>}[dd]^{\txt{No}} 
& 
\algbox[\lrat]{\center Calculate $D_0$.\\ 
 Is $D_0\equiv0$?}
\ar@{=>}[ur]^{\txt{No}}   \ar@{=>}[r]^{\txt{Yes}}
& 
\algbox[\rrat]{\noindent{\bf b1} \hfil
                $F$ is given\ \ \hfill\\ 
            \hbox{\ \ }  by Proposition \ref{R=const}.}
\\
&&
\algbox[\rrat]{\center For given $g$ and $A$\\ $F$ does not exist.}
\\
\algbox[\lrat]{\center Calculate $\varphi_2$.\\ Is $\varphi_2\equiv0$?} 
\ar@{=>}[r]^{\txt{No}}\ar@{=>}[d]^{\txt{Yes}}
&
\algbox[\lrat]{\center 
    Calculate $\calg_2,\calg_3$.\\
   Are $\calg_2\equiv0\equiv\calg_3$?} 
  \ar@{=>}[ur]^{\txt{No}}
  \ar@{=>}[r]^{\txt{Yes}}
&
\algbox[\rrat]{\center $F$ is given\\ by Theorem \ref{th:main1}.}
\\
\algbox[\lrat]{\center  Calculate $\cald_x,\cald_y$.\\
{\bf b2} \hfill  Is $\cald_x\equiv 0 \equiv \cald_y$?\ \ \ \hfill}
\ar@{=>}[r]^{\txt{No}}
\ar@{=>}[d]^{\txt{Yes}}
&  
\algbox[\lrat]{\center For given $g$ and $A$\\ $F$ does not exist.}
&
\\
\algbox[\lrat]{\center Calculate $\varphi^*_2$.\\ Is $\varphi^*_2\equiv0$?} 
\ar@{=>}[r]^{\txt{No}}\ar@{=>}[d]^{\txt{Yes}}
&
\algbox[\lrat]{\center 
    Calculate $\calg^*_2,\calg^*_3$.\\
   Are $\calg^*_2\equiv0\equiv\calg^*_3$?} 
  \ar@{=>}[dr]^{\txt{No}}
  \ar@{=>}[r]^{\txt{Yes}}
&
\algbox[\rrat]{\center $F$ is given\\ by Theorem \ref{th:main2}.}
\\
\algbox[\lrat]{\center $g$ admits a  Killing\\
{\bf b3}\hfill vector field. \hfill } 
\ar@{=>}[d]
&
&
\algbox[\rrat]{\center For given $g$ and $A$\\ $F$ does not exist.}
\\
\algbox[\lrat]{\center Calculate $\cald^*_x,\cald^*_y$.\\
{\bf b4} \hfill Is $\cald^*_x\equiv \cald^*_y\equiv0$?\ \ \ \hfill}
\ar@{=>}[r]^{\txt{No}}
\ar@{=>}[d]^{\txt{Yes}}
&  
\algbox[\lrat]{\center For given $g$ and $A$\\ $F$ does not exist.}
&
\\
\algbox[\lrat]{\center $F$ is given by \\ Proposition \ref{A-kill}.}
&
}
\]}
\end{minipage}
\pagebreak[2]

Let us comment some boxes in the flowchart of the algorithm.  
Recall that \slsf{$A$ is compatible with $g$} if $A$ is the
Birkhoff-Kolokoltsov $3$-codifferentials associated to some cubic integral $F$
of the metric $g$.

\begin{itemize}
\item[\bf(b1)] If $R$ is constant, then every cubic integral is of the form
 $F=\alpha^{ijk}L_iL_jL_k$ with constants $\alpha^{ijk}\in\rr$ where $L_i$, $i=1,2,3$ are
 three linear (polynomial of degree 1) invariants which are linearly
 independent. Thus the cubic integrals  form a linear space of dimension $10$, see
 \cite{Kr}. On the other hand, the space of holomorphic $3$-codifferentials
 $A$   has infinite dimension.  Proposition \ref{R=const} shows that the
 equation $D_0=0$ is necessary and sufficient condition for compatibility of
 $A$ with $g$ and, if it is fulfilled, the space of cubic integrals depends on
 $3$ real parameters. More precisely, the a generic cubic integral has
  the form $\Re(A)+b^ip_i\cdot H$ where $\vec b=(b^i)$ is a vector field with given
  values $\vec{b}(P)$ and $\mathsf{rot}\,\vec{b}(P)$ at a given point $P\in
  S$. The construction of the
  vector field $\vec{b}=(b^i)$ reduces to solution of  certain linear systems
  of ODEs.

\item[{\bf(b2)}] It is sufficient to calculate only one of the
  expressions $\cald_x,\cald_y$, namely one such the corresponding
  $\tfrac{\partial\varphi_0}{\partial x}$, $\tfrac{\partial\varphi_0}{\partial y}$ is not zero, see Theorem
  \ref{th:main2}.  
 
\item[\bf(b3)] By a classical result of Darboux and Eisenhart (see discussion
 in \S \ref{previ}) the vanishing of both $\varphi_2$ and $\varphi_2^*$
 is equivalent to the local existence of a Killing vector field  which we
  denote by $\xi$. The latter is equivalent to the existence of a linear
 integral  which we denote by $L$. Since Lie derivative of the curvature
 $R$ along $\xi$ must vanish, we have that $\xi$ is  proportional to the
 Hamiltonian vector field of $R$, $\xi=f(x,y)\cdot\vec X_R$. Such $f(x,y)$ is unique
 up to a constant factor, and can be found from the  PDE-system
  $\scrl_{f\vec X_R}g=0$ on the unknown function $f$. Note that the PDE-system
  $\scrl_{f\vec X_R}g=0$ is Frobenius and, in view of conditions $\varphi_2=0$ and
  $\varphi_2^*=0$, is in involution, i.e., is essentially a system of ODE.

\item[{\bf(b4)}] The expressions $\cald^*_x,\cald^*_y$ and the part of the
 algorithm starting from the box {\bf(b3)} is treated in §\,\ref{Killing}.

\end{itemize}

In the pseudo-Riemannian this algorithm basically works, however there exists
an important difference. Namely, the condition $|\nabla R|^2_g=0$ does not imply
that $R=\const$.  As a consequence we obtain a new special class of
pseudo-Riemannian metrics whose counterpart does not exists in the Riemannian
case.  They are not metrics of constant curvature, they have $\varphi_2\equiv0\equiv\varphi_2^*$,
however, as we prove in Proposition \ref{qub-phi1}, the metrics from this
class admit no cubic integrals, and therefore no Killing vector fields.  We
refer to Section \ref{sign=1,1} for further details and description of an
additional step in the algorithm.

We emphasise here that the algorithm works in any (not necessary isothermal)
coordinate system.  In particular, Lemma \ref{A-form} gives a way how to
verify whether the data $(g,A)$ have sense, i.e., whether \emph{ real}
(3,0)-tensor $\hat A=\Re(A)$ corresponds to the real part of a holomorphic
3-codifferential.

\subsection{ Motivation and history} \label{intr2}
  
The importance of polynomial integrals for studying the geodesic flow was
recognised long time ago. Indeed, it was Jacobi who had realised that the
geodesic flow on the ellipsoid admits an ``extra'' quadratic integral. This
allowed Jacobi to integrate the geodesic flow on the ellipsoid.

Polynomial integrals for the geodesic flow are interesting for physics and for
differential geometry. The interest from physics is due to the following
observation (called the Maupertuis principle, or coupling constant transform):
for sufficiently large energy levels $h$ the restriction of the Hamiltonian
system on $T^*M$ with the Hamiltonian $ H:=\half|\vec p|^2 + V $ (where
$\half|\vec p|^2$ is again the kinetic energy corresponding to $g$, and
$V:{S}\to \mathbb{R}$ is the potential) to the energy surface %
$\{(x,\vec p)\in T^*M\mid\;H(x,\vec p)= h\}$ has the same unparametrised orbits as
the geodesic flow of the metric $g_h:= (h- V)g$.  Moreover, if the Hamiltonian
system is polynomially-integrable, then the geodesic flow is also
polynomially-integrable.  More precisely, if a function
\begin{equation} \label{eq:f1} 
F(x,y;p_x,p_y)= 
\underbrace{a_0 p_x^3 + a_1  p_x^2p_y + a_2 p_x p_y^2 + a_3 p_y^3
}_{F_3(x,y;p_x,p_y)} + 
\underbrace{c_0 p_x   +c_1 p_y}_{F_1(x,y;p_x,p_y)}
\end{equation}  
is an integral of the Hamiltonian system with the Hamiltonian
$\half|\vec{p}|^2+V$, then the (homogeneous cubic in momenta) function
\begin{equation} \label{eq:mp} 
F_h:= F_3 + H_h \cdot  F_1 
\end{equation} 
is an integral for the Hamiltonian system with the Hamiltonian
$H_h:=\frac{1}{2(h-V)}|\vec p|^2$, i.e., for the geodesic flow of the metric
$g_h:=(h-V)g$.
  
Note that the assumption that the integral \eqref{eq:f1} does not contain the
terms $F_0$ and $F_2$ of degree $0$ and $2$ in momenta is a natural one: by
Whittaker \cite{Wh}, if a function of the form $F_0+ F_1 + F_2 +F_3$ is an
integral for the Hamiltonian system with the Hamiltonian
$H:=\half|\vec{p}|^2+V$, then $F_1+F_3$ is also an integral.
     
\smallskip%
Metrics admitting integrals polynomial in velocities are also interesting for
geometry, since the existence of such integrals implies interesting geometric
properties of the metric. Moreover, study of polynomial integrals helps to
solve pure geometrical problems. For example, the existence of the integral of
degree 1 is equivalent to the existence of a Killing vector field, and implies
the existence of the coordinate system such that the component of the metric
does not depend on one of the coordinates \cite{Da}.  The existence of an
integral of degree 2 is equivalent to the existence of another metric having
the same geodesic with the initial one \cite{Ma-To1,Ma-To2}.  The last
observation was recently used \cite{Br-Ma-Ma,Ma} in the solution of two
geometric problems explicitly formulated by Sophus Lie in 1882.

It is not clear whether the existence of the integral of degree 3 has such a
clear geometric meaning, but still, as we mentioned above, it allows one to
construct a holomophic 3-codifferential with interesting properties, and, as
it was shown in \cite{Du-Ma-To}, a volume-preserving vector field on the
surface.

\subsection{Previous results in this field} \label{previ} 

The most classical result in this direction is in folklore attributed to
Bonnet, and appeared in the books of Darboux \cite{Da} and of Eisenhart
\cite{Ei}.
 
Given a metric $g =(g_{ij}(x,y))$ on a surface, let us consider the functions
$\varphi_0,\varphi_1,\varphi_1^*,\allowbreak\varphi_2,\varphi_2^*$ from \S \ref{intr1}: $\varphi_0:=R$ is the
Gauss curvature, $\varphi_1:=\half|\nabla R|^2_g$ is half-square of the gradient of $R$,
$\varphi_1^*:=\Delta_gR=\Delta_g\varphi_0$ the metric Laplacian of $R=\varphi_0$, and $\varphi_2:=\{\varphi_0,\varphi_1\}_g$,
$\varphi^*_2:=\{\varphi_0,\varphi^*_1\}$ are the Poisson brackets.  

One can show (see Darboux \cite[§§688, 689]{Da} and Eisenhart
\cite[pp.\;323-325]{Ei}) that the metric admits a nontrivial integral linear
in momenta (locally, in a neighbourhood of almost every point) if and only if
the invariants $\varphi_2$, $\varphi_2^*$ vanish.

Though Darboux and Eisenhart proved this result for Riemannian case only,
their proof can be generalised to the pseudo-Riemannian case with the
following reformulation of the assertion: {\it if the equations $\varphi_2=0$ and
 $\varphi^*_2=0$ are satisfied and $\half|\nabla R|^2_g=\varphi_1$ is non-zero, then $g$ admits a
 nontrivial linear integral.}  We study the class of pseudo-Riemannian metrics
with the condition $|\nabla R|^2_g\equiv0$ in Section \ref{sign=1,1}, and prove
 that they admit a linear integral iff $R=\const$.

We see that we have an algorithmic  way to decide
whether a given metric admits a linear integral of degree: one need to
calculate the invariants $\varphi_2$, $\varphi_2^*$ which are algebraic expressions in the
components of the metric and their derivatives up to order 5, and compare them
with zero.

It was a long standing problem to find a similar algorithmic way to decide
whether a given metric admits an integral quadratic in momenta. The first
attempts are due to Roger Liouville%
\footnote{Roger Liouville was a younger relative of the more famous Joseph
 Liouville and attended his lectures at the Ecole Polytechnique}
\cite{Li} and Koenigs \cite{Koe}.  They understood that the PDE-system for the
coefficients of the integral is linear and of finite type, implying that that
there must exist an algorithmic way to decide whether a given metric admits an
integral quadratic in momenta.  Unfortunately, the calculation were too
complicated to be fulfilled in that time.  An interesting approach to this
problem is due to Sulikovski \cite{Su}: he invented a trick which allowed to
simplify the calculation.  Finally, the problem was solved in recent
independent works by Kruglikov \cite{Kr} and Bryant et al
\cite{Br-Du-Ea}.

Geodesic flows admitting cubic integrals is a much more difficult subject than
those admitting linear or quadratic integrals.  Indeed, metric admitting
linear and quadratic integrals are completely described: the local description
is done by classics \cite[§§540--596]{Da}, see \cite{Bo-Ma-Pu1,Bo-Ma-Pu2} for the
 discussion of the pseudo-Riemannian case, and the global (= if the surface
is closed) is due to Kolokoltsov \cite{Kol} and Kiyohara \cite{Ki1}, see the
survey \cite{Bo-Ma-Fo}. The metrics admitting linear integrals are
parameterised by one function of one variable, the metrics admitting quadratic
integrals are parameterised by two functions of one variable.  No such
description  is known in the cubic case (though from general theory of
PDE it follows that locally, in the analytic formal category, the metrics
admitting cubic integrals are parametrised by 3 functions of one variable).

Though the complete description of metrics admitting cubic integrals is not
known yet, there exists a lot of local examples  (the first local
 examples appeared already in Darboux \cite[§619]{Da}).  Some of most
interesting examples come from physics, see the survey \cite{Hi}, or are
motivated by physics \cite{Ha,Ye,Ve-Ts}.

It appears to be difficult to construct metrics admitting cubic integrals on
closed surfaces. As it was shown by Kozlov \cite{Koz} and Kolokoltsov
\cite{Kol}, the surfaces of genus $\geq 2$ do not admit (nontrivial) polynomial
integrals; we repeat the argument of Kolokoltsov in  the proof of 
Corollary \ref{g-geq-2}.  Of course, the metrics admitting linear integral
also admit cubic integrals, since the cube of the linear integral is a cubic
integral. Besides these half-trivial examples, the list of known examples is
very short: there is a series of explicit examples coming from Goryachev case
of rigid body motion and its generalisation, a series of recently found
examples due to \cite{Du-Ma}, and proof of the existence of other examples in
\cite{Ki2}, \cite{Se}.

In particular, in all known examples the surface is the sphere (i.e., there is
no known example of metric on the torus admitting cubic integral and do not
admitting linear integral; moreover, Bolsinov et al conjectured
\cite{Bo-Ko-Fo} that such examples can not exist).

The construction of curvature-type invariants whose vanishing is equivalent to
the existence of a cubic integral is discussed by Kruglikov in \cite[§12]{Kr}.
In particular, he has shown that the dimension of the space of cubic integrals
is $\leq 10$, moreover, the metrics admitting 10-dimensional space of cubic
integrals have constant curvature. It was conjectured by Kruglikov in \cite{Kr}
that the next largest value of the dimension of the space of cubic integrals
is $4$ and that in this case $g$ is Darboux superintegrable
(the definition see in \cite{K-K-M-W} or
in \cite[\S 2.2.4]{Br-Ma-Ma}). If the Kruglikov conjecture is true,
\cite[Theorem 2]{Kr} gives curvature-type invariants whose vanishing is
equivalent to the existence of the 4-dimensional space of cubic integrals.

In the present paper, we answer the question whether the pair (metric $g$,
holomorphic 3-codifferential $A$) is compatible in the sense there exists a
cubic integral for $g$ whose Birkhoff-Kolokoltsov form is $A$. On the level on
PDE, we reduce the number of unknown functions, which of course make the
problem easier. This is the reason that our answer (for cubic integrals)
is more simpler than the answers of \cite{Kr} and \cite{Br-Du-Ea} (for
quadratic integrals).  From other side, the assumption that $A$ is given is
natural from the viewpoint of classical mechanics.  Indeed, all integrals
$F_h$ corresponding to different values of $h$ (see \S \ref{intr2}) have
the same holomorphic 3-codifferential $A$.  Moreover, on closed surfaces,
the space of holomorphic 3-codifferentials is finite-dimensional. That means
that, in the investigation of cubic integrals on closed surfaces, one can view
our compatibility conditions as the conditions on the metric only depending on
finitely many parameters.

\section{Complex calculus and prolongation-projection method.}
\label{sec-calcul} 

\subsection{Homogeneous polynomial integrals.} \label{F-homo} Let $g$ be a
metric on a surface $S$ and let $H=\half |p|^2_g$ be a Hamiltonian given by
the kinetic energy and let $F(x,p)$ be its integral which is polynomial of
degree $\leq d$ in momenta.  Then $F(x,p)=\sum_{j=0}^dF_j(x,p)$ where each component
$F_j$ is homogeneous of degree $j$ in momenta.  Direct computation shows that each Poisson
bracket $\{H,F_j\}$ is homogeneous of degree $j+1$ in momenta. The uniqueness of
the decomposition of a polynomial into homogeneous components implies that
the Poisson bracket $\{H,F\}$ vanishes if and only if each bracket $\{H,F_j\}$
vanishes. By this reason, considering polynomial integrals of purely kinetic
Hamiltonians $H=\half |p|^2_g$ we can restrict ourselves to homogeneous
polynomial integrals $F$.

\subsection{Birkhoff-Kolokoltsov codifferential.}\label{cplx-calc}
Every metric $g$ on a surface $S$ admits locally so-called \slsf{isothermic
 coordinates} $(x,y)$ in which the metric has the form
$g=\lambda(x,y)(dx^2+dy^2)$. If we orient $S$ by means of these coordinates, then
the function $z:=x+\cplxi y$ will be a complex coordinate for the \slsf{induced
 complex structure} on $S$, characterised by the following property: The
multiplication of $z$ by $\cplxi$ is the rotation of the chart $(x,y)$ by
$90^\circ$.  The existence of isothermic coordinates was discovered by
B.~Riemann. A detailed proof under very weak assumption on the metric (e.g.,
$g$ is merely continuous) can be found in \cite{Beg}.

In forthcoming calculus we use such a complex coordinate $z$ and its conjugate
$\bz$ instead of real coordinates $(x,y)$, since most formulas become
simpler and more compact. The corresponding dual  coordinates are  
$p:= \frac12(p_x-\cplxi p_y)$ and $\barp := \frac12(p_x+ \cplxi p_y)$, the
function $H$ is expressed as $H=\frac{2p\barp}{\lambda}$, and the canonical
symplectic form in the cotangent bundle $T^*S$ (``phase space'') is expressed
as $\omega_{can}=dp\land dz+d\barp\land d\bz$.

Now assume that $F=F(z,\bz, p, \barp)$ is a (local) integral of $H$ cubic in
momenta. Then $F= \Re( a\cdot p^3 + b \cdot p^2\barp)$ where $\Re$ denotes the real part
of a complex expression and $a=a(z,\bz)$, $b=b(z,\bz)$ are some smooth
(local) \emph{complex-valued} functions on $S$. An explicit calculation (see
\cite{Du-Ma-To}) yields:
\begin{equation}\label{braFH}
\textstyle
\{F,H\} = \frac{2}{\lambda^2}\Re\big(
p^4{\cdot} a_\bz + p^3\barp{\cdot}( b\lambda_\bz +3a\lambda_z+ \lambda b_\bz+\lambda a_z)+
p^2\barp^2 {\cdot} (2b\lambda_z + \lambda b_z)
\big),
\end{equation}
where $(\cdot)_z$ and $(\cdot)_\bz$ denote the complex derivatives: 
\[\textstyle
f_z = \frac12\left(\frac{\partial f}{\partial x} - \cplxi \frac{\partial f}{\partial y}\right)
\qquad
f_\bz = \frac12\left(\frac{\partial f}{\partial x} + \cplxi \frac{\partial f}{\partial y}\right)
\]

Since the term $p^4\cdot a_\bz$ must vanish, we obtain immediately the following

\begin{lem}\label{A-hol} In the situation above, if $F=\Re(a\cdot p^3+b\cdot p^2\barp)$ is
 an  integral, then the coefficient  $a$ is a 
 holomorphic.
\end{lem}

Rewrite any function $F(x^1,x^2;p_1,p_2)$ cubic in momenta as a polynomial
$F=F^{ijk}p_ip_jp_k$ in $p_i$. Then, under change of coordinates, its
coefficients $F^{ijk}(x^1,x^2)$ transform as components of a symmetric
$(3,0)$-tensor. The same holds for the ``$a$-component'' of $F$: Under
\emph{complex} (\,$\Leftrightarrow$ holomorphic) coordinate change $z\mapsto z'$ it transforms as
$a'= a\cdot\left(\frac{\partial z'}{\partial z}\right)^3$.  Indeed, the momenta $p$ and
 $\bar p$ are transformed as $p= p' \tfrac{dz'}{dz}$, $\bar p= \bar p'
 \tfrac{d\bar z'}{d\bar z}$ implying
\[
F = \Re(a\cdot p^3+b\cdot p^2\barp)= \Re\Big(\underbrace{a \cdot\left( \tfrac{dz'}{d
    z}\right)^3 }_{a'} \cdot {p'}^3+b \cdot \left( \tfrac{dz'}{dz}\right)^2 \,
 \tfrac{d\bar z'}{d\bar z} \cdot {p'}^2\barp'\Big).
\]
This means that
$A:=a(z)\,\big(\frac{\partial}{\partial z}\big)^{ {\otimes3}}$ is a well-defined section of the bundle
$(T^{(1,0)}S)^{\otimes3}$ independent of the choice of a local holomorphic
coordinate.  By Lemma \ref{A-hol}, the section is holomorphic. 

\begin{defi}\label{3-codiff} \rm For any function $F(x^1,x^2;p_1,p_2)$ cubic in
 momenta the complex tensor $A:=a\cdot \big(\frac{\partial}{\partial z}\big)^3$ is called the
 \slsf{Birkhoff-Kolokoltsov 3-codifferential} associated with $F$.  We shall
 usually denote by $a=a(x,y)=a(z,\bz)$ its complex-valued coefficient (wrt.\
 some complex coordinate $z$), and denote by $\hat A:=\Re(A)$ its real part. The
 latter is a symmetric $(3,0)$-tensor.
\end{defi} 

Since every non-vanishing holomorphic function $f(z)$ in one variable has
isolated zeroes, the same is true for any holomorphic 3-codifferential. Notice
that if $a(z)$ is holomorphic and non-vanishing, then locally near a given
point $p$ there exists the root $\sqrt[3]{a(z)}$ which is again a
non-vanishing holomorphic function. It follows that making the change to the
new complex coordinate $z'$ given by $z':=\int\frac{dz}{\sqrt[3]{a(z)}}$, the
component $a(z)$ transforms into the constant $a'(z')\equiv1$.

\smallskip %
Another applications of holomorphicity of $A$ are as follows:

\begin{corol}[\cite{Kol}] \label{g-geq-2}   \sli
 No Riemannian metric on a closed oriented surface $S$ of genus $g\geq2$ admits a
 global non-trivial polynomial integral $F$.

\slii Let $S$ be the $2$-torus $T^2$, $g$ a metric, $z$ a \emph{linear}
holomorphic coordinate, and  $F=\Re(a\cdot p^3+b\cdot p^2\barp)$ a cubic integral. Then
$a$ is constant. 
\end{corol} 

\proof \sli As is is shown in §\,\ref{F-homo}, we may assume that our
invariant $F$ is homogeneous of degree $k$. Evidently one
can take $F$ in the form $F(z,\bz,p,\bar p)= \Re( a\cdot p^k) + p\barp\cdot B$ where
$B(z,\bz,p,\bar p)$ is homogeneous in momenta of degree $k-2$. An explicit
computation shows that similarly to the cubic case %
$\{H,F\}= \frac{2}{\lambda^2} \Re(a_\bz\cdot p^{k+1})+p\barp\cdot C$ for some function
$C(z,\bz,p,\barp)$ which is homogeneous polynomial of degree $k-1$ in
momenta. Thus $a_\bz$ must vanish which means that $a(z)$ is holomorphic as in
the cubic case. Similarly to the cubic case, %
$A:=a\cdot \big(\frac{\partial}{\partial z}\big)^k$ is a well-defined \emph{holomorphic} section
of the bundle $(T^{(1,0)}S)^{\otimes k}$.

The classical complex geometry (see e.g., \cite{Gr-Ha}) says that for any
non-zero \emph{meromorphic} section $A$ of the bundle $(T^{(1,0)}S)^{\otimes k}$
the difference $n(A)-p(A)$ between the numbers of zeroes and poles of
$A$ is $k(2-2g)$. Since $A$ is holomorphic in our case, there are no
poles and $k(2-2g)$ must be non-negative. Consequently, in the case $g\geq2$
$A$ must vanish identically. It hollows that $\frac{1}{\lambda}B=F\cdot H\inv$ is a
global homogeneous integral for $H$ of degree $k-2$. Now we can apply
induction in $k$.

\medskip %
\slii A local holomorphic coordinate on a torus $T^2$ is {\it linear}
if it comes from the standard complex coordinate $z$ on the complex plane
$\cc$ under some conformal isomorphism $T^2\cong\cc/\Lambda$ for some 2-lattice $\Lambda\subset\cc$.
An equivalent condition is that $\frac{\partial}{\partial z}$ is a globally defined
holomorphic vector field. Such a field on a torus
is non-vanishing, and therefore $a$ is a globally defined holomorphic function
on a closed Riemann surface. Such a function must be constant by Liouville's
theorem.
\qed

\subsection{Prolongation-projection method. } \label{prol1}
The ``naive'' idea behind this method is very simple. Suppose we have a system
$\cals=\{\cale_i\}$, of PDEs which is \slsf{overdetermined,} i.e., the number of
equations is bigger than the number of unknown functions.  In our case, we
have two equations on one unknown function. We \slsf{prolong } the system,
i.e., we differentiate the equations with respect to all variables (in our
case, $x$ and $y$) and add the results to the system $\cals$. In other words,
we consider the systems
\[\textstyle
\cals^{(1)}:=
\{\cale_i \}\bigcup \{\tfrac{\partial}{\partial x}\cale_i \} \bigcup
\{\tfrac{\partial}{\partial y} \cale_i \}
 \ \ \textrm{ --- first prolongation,
}
\]
\[\textstyle
\cals^{(2)}:=
\{\cale_i \}\bigcup \{\tfrac{\partial}{\partial x}\cale_i \} \bigcup \{\tfrac{\partial}{\partial y} \cale_i \}
\bigcup \{\tfrac{\partial^2}{\partial x^2}\cale_i \}\bigcup \{\tfrac{\partial^2}{\partial x\partial y}\cale_i \}\bigcup 
\{\tfrac{\partial^2}{\partial y^2}\cale_i \}
 \ \ \textrm{ --- second prolongation,
}
\]
and so on. Obviously, the operation ``prolongation'' does not change
the set of sufficiently smooth solutions of the system.  Indeed, every
(sufficiently smooth) solution of $\cals$ is also a solution of the
prolonged systems.  From the other side, every solution of the prolonged
system is evidently a solution of $\cals$, because $\cals$ is a
part of the prolonged system.

If the system is overdetermined, then generically the number of new equations
grows faster than the number of derivatives of the unknown functions. Then, 
in the generic case, after say $l$ prolongations one can
resolve the prolonged system $\cals^{(l)}$ with respect to the highest
derivatives. In this case, the system $\cals$ is called  the systems of
 \slsf{finite type.}  In our case, the system is indeed of finite type, and the
number of needed prolongations is $l=1$.

After this point, we can start the prolongation-\slsf{projection} procedure.
The key observation here is as follows: The successive prolongation
$\cals^{(l+1)}$ can be still resolved with respect to the highest derivatives,
the number of the equations in $\cals^{(l+1)}$ will be generically greater
than the number of highest derivatives, and with the help of algebraic
manipulations we can obtain \emph{new} equations of \emph{lower} degree.  This
algebraic manipulation --- expressing of higher derivatives via lower ones
using part of equations and substituting these expression in remaining
equations --- is called \slsf{projection} procedure.  We can repeat the
prolongation-projection many times (clearly, is sufficient to differentiate
only \emph{newly obtained} equations).
 
In the generic case, after finite number of prolongation-projections,  we obtain
a system in which the number of \emph{algebraically independent} equations is
greater than the number of partial derivatives of unknown functions (we consider the unknown functions  as partial derivatives of order $0$). Such an algebraic system is
inconsistent. In this case the initial system $\cals$ has no solution. It may  also 
happen that certain  prolongation $\cals^{(k)}$, considered as a system of algebraic equations    on the derivatives of the unknown functions,  is algebraically inconsistent, even if the number of algebraically independent is less than the number of partial derivatives.

Thus, the existence of a solution of an overdetermined system $\cals$ of
finite type implies that after certain number of prolongations and
prolongation-projections we come to the point  where prolongation-projections
do not produce essentially new equations. The latter means that the equations 
we obtain are algebraic corollaries of the equations we already have. Moreover,  
the system of algebraic equations on the partial derivatives of
unknown functions is consistent.  In this case, the system is called 
 \slsf{involutive}, or \slsf{in involution}.  One can show that (under certain
additional assumptions which are fulfilled in our case) involutive systems can
be solved locally. Moreover, one can find all solutions with the help of
solving  of certain ODEs and of algebraic operations.

As we noted above, a generic overdetermined system   of finite type is 
inconsistent and can not be put to be in involution by
prolongation-projections.  One obtains a system in involution after certain
number of prolongations and prolongation-projections, if the coefficients of
the initial system satisfies certain partial differential equations (called
\slsf{integrability conditions}).  These PDE on the coefficients are
 equivalent to two conditions:
\begin{itemize} 
\item[(1)] we obtain no new equations after certain prolongation-projection.
\item[(2)] All systems  we obtain by prolongations and prolongation-projections, considered as  algebraic systems on the partial 
 derivatives of the unknown functions,  are  algebraic consistent.
\end{itemize}

In our paper, we find explicitly all integrability conditions for the system
that corresponds that $g$ and $A$ are compatible: they are $\calg_2,\calg_3$,
$\calg_2^*$,$\calg_3^*$, $\cald_x,\cald_y$, $D_0$ from Theorems
\ref{th:main1}, \ref{th:main2} and Propositions \ref{R=const}, \ref{A-kill}.

The prolongation-projection method is a very powerful method for finding
solutions of overdetermined systems of PDE of finite type. Moreover, in our
case the initial system of PDE is an inhomogeneous linear system, implying that
standard difficulties due to impossibility of solving explicitly systems of
algebraic equations do not appear, so the method can be applied
algorithmically.  Actually, the only difficulty in applying the
prolongation-projection method for linear overdetermined systems of PDE of
finite type is the calculational one.  In many cases, even modern computer
algebra programs are not powerful enough for using prolongation-projection.

We overcome this difficulty with the help of ``advanced complex calculus''
introduced in \S \ref{calcul}: roughly speaking, this calculus is a
collection of quite nontrivial tricks that allows us to `hide' the covariant
derivatives of the objects, and, therefore, makes all formulas very compact,
so we could do prolongation-projections ``by hands'', without using computer
algebra. Moreover, all differential equations on the coefficient of the metric
and of the Birkhoff-Kolokoltsov form are automatically expressed in the
invariant form, so one can calculate them in an arbitrary coordinate system.

\smallskip%
\state Notation. Applying the prolongation-projection method we use
 the following terminology.  Let $\cale$ be a (differential) equation or
 expression and $\cals=\{\cale_i\}$ a system (ie., a set) of such equation. Then
 $\cale$ is an \slsf{algebraic consequence} of $\cals$ if it can be deduced
 from $\cals$ by means of algebraic manipulations, and a \slsf{differential
  consequence} if $\cale$ is an algebraic consequence of $\cale$ and partial
 derivatives of equations from $\cals$. In this paper we work with
 \emph{linear} PDEs and make linear-algebraic manipulation with equations. In
 particular, we have \slsf{linear-algebraic consequences} of linear systems
 of equations, and differential
 consequences can be obtained applying \emph{linear} differential operators.

\smallskip%
\state Remark. The name ``prolongation-projection'' is natural in the context of
 jet bundle geometry, we recommend the textbook \cite{Kr-Ly-Vi} as a
 source. Equivalent approach is called Cartan, or Cartan-Kähler theory, we
 recommend the textbooks \cite{B-C-G-G-G,Iv-La} for more details.  The system
 of PDE we consider in our paper is easier than the generic systems treated in
 \cite{Kr-Ly-Vi,B-C-G-G-G,Iv-La} (because they are linear and of finite type),
 that's why we can achieve our goals without going too far into the geometric
 theory of PDE. Actually, our intention is to make the paper understandable
 for a standard mathematician working in mathematical physics; that's why we avoided the
  terminology   of the geometric
 theory of PDE.

\subsection{Advanced complex calculus.}\label{calcul}
Our next goal is to develop certain calculus which will allow us to simplify
computation and express results and formulas in more compact form. 

Let $(S,g)$ be a Riemannian surface with a fixed orientation and
$z=x+\cplxi y$ a local complex coordinate. Let $T^*S^\cc$ be the complexified
cotangent bundle, its sections are complex $1$-forms on $S$. Then sections %
$dz=dx+\cplxi dy$ and $d\bz=dx-\cplxi dy$ form a local frame of
$T^*S^\cc$. Consider the decomposition $T^*S^\cc=\Omega^{(1,0)}S\oplus\Omega^{(0,1)}S$ in
which the summands are generated by the forms $dz$ and respectively $d\bz$.
The formula for transformations of the forms $dz$ and $d\bz$ under
\emph{holomorphic} coordinate changes shows that the
decomposition $T^*S^\cc=\Omega^{(1,0)}S\oplus\Omega^{(0,1)}S$ is independent of the choice of
a local complex coordinate. A similar decomposition $TS^\cc=T^{(1,0)}S\oplus
T^{(0,1)}S$ %
for the complexified tangent bundle is obtained using local vector fields
$\frac{\partial}{\partial z}:=\frac12(\frac{\partial}{\partial x} - \cplxi \frac{\partial}{\partial y})$ and $\frac{\partial}{\partial
 \bz}:=\frac12(\frac{\partial}{\partial x} + \cplxi \frac{\partial}{\partial y})$.

We refer to \cite{Gr-Ha,For} for basic properties of these
decompositions. Main properties of the bundles  $T^{(1,0)}S$, $T^{(0,1)}S$, $\Omega^{(1,0)}S$,
$\Omega^{(0,1)}S$ which we shall use are as follows:
\begin{enumerate}
\item The bundles $\Omega^{(1,0)}S$ are $\Omega^{(0,1)}S$ are \emph{complex dual
  line bundles} to $T^{(1,0)}S$ and $T^{(0,1)}S$, respectively. 
\suspend{enumerate}
Using this fact, we define complex line bundles $T^{[p,q]}S$ and $\Omega^{[p,q]}S$
with \emph{integer} $p,q$
setting 
\[
T^{[p,q]}S := \big(T^{(1,0)}S\big)^{\otimes p} \otimes \big(T^{(0,1)}S\big)^{\otimes q},
\qquad
\Omega^{[p,q]}S := \big(\Omega^{(1,0)}S\big)^{\otimes p} \otimes \big(\Omega^{(0,1)}S\big)^{\otimes q}.
\]
These bundles should be not confused with bundles of tensors of type $(p,q)$
used in differential geometry, and with bundles of $(p,q)$-forms used in
complex analysis and geometry. In order to emphasise the difference, we use
brackets $[p,q]$ instead of parentheses $(p,q)$.

\resume{enumerate} 
\item[(1')] From the definition and property (1) we obtain immediately 
\[
\Omega^{[p,q]}S=T^{[-p,-q]}S
\quad\text{\ and \ }\quad
T^{[p,q]}S
 \otimes T^{[p',q']}S =T^{[p+p',q+q']}S.
\]
The bundle $T^{[0,0]}=\Omega^{[0,0]}$ is trivial and its sections are usual
complex-valued functions.

\item\label{pr2}
 The bundle $(T^*S)^{(\otimes k)}\otimes\cc$ (respectively $(TS)^{(\otimes k)}\otimes\cc$) splits
 into the sum of line bundles each isomorphic to some $\Omega^{[p,q]}S$
 (respectively $T^{[p,q]}S$) with $p,q\geq0$ and $p+q=k$. A similar splitting
 holds for a general complex tensor bundle $(TS)^{(\otimes k)}\otimes(T^*S)^{(\otimes l)}\otimes\cc$.
 The components of $(T^*S)^{(\otimes k)}\otimes\cc$ can be indexed by the order of factors
 $\Omega^{(1,0)}S$, $\Omega^{(0,1)}S$ in the product.

 Notice that the subbundle of \emph{symmetric tensors} of $(T^*S)^{(\otimes k)}\otimes\cc$
 is isomorphic to the sum $\oplus_{p=0}^k\Omega^{[p,k-p]}S$, so that each $\Omega^{[p,k-p]}S$
 appears exactly once. In the special case $k=2$ the power $(T^*S)^{(\otimes2)}\otimes\cc$
 has as summands one bundle $\Omega^{[2,0]}$, one bundle $\Omega^{[0,2]}$, and two
 bundles $\Omega^{(1,0)}S\otimes\Omega^{(0,1)}S$, $\Omega^{(0,1)}S\otimes\Omega^{(1,0)}S$ each isomorphic to
 $\Omega^{[1,1]}$.

Moreover, for \emph{positive} $i\leq\min(p,k),j\leq\min(q,l)$, the natural
 isomorphism
\begin{equation}\label{trace}
T^{[p,q]}S\otimes\Omega^{[k,l]}S \xrar{\;\cong\;}  T^{[p-i,q-j]}S\otimes\Omega^{[k-i,l-j]}S
\end{equation}
coincide with the $(i+j)$-fold index contraction in the space of symmetric
complex-valued $(p+q,k+l)$-tensors.

If $z=x+\cplxi y$ is a local complex coordinate on $S$, then
\[
dz\otimes d\bz= (dx^2+dy^2)+\cplxi(dx\otimes dy-dy\otimes dx)=(dx^2+dy^2)+\cplxi dx\land dy
\]
and $d\bz\otimes dz=(dx^2+dy^2)-\cplxi dx\land dy$. Thus the symmetric part of the sum
$\Omega^{(1,0)}S\otimes\Omega^{(0,1)}S\oplus\Omega^{(0,1)}S\otimes\Omega^{(1,0)}S$ is spanned by
$dz\bullet d\bz:=\frac12(dz\otimes d\bz+d\bz\otimes dz)=dx^2+dy^2$ and its antisymmetric part
spanned by $dz\land d\bz=2\cplxi dx\land dy$. This gives another decomposition of sum
$\Omega^{(1,0)}S\otimes\Omega^{(0,1)}S\oplus\Omega^{(0,1)}S\otimes\Omega^{(1,0)}S$ into two bundles each isomorphic
to $\Omega^{[1,1]}$. The anti-symmetric part is the bundle $\Lambda^2S\otimes\cc$ of complex
$2$-forms on $S$ and the power $(T^*S)^{(\otimes2)}\otimes\cc$ splits into the sum of
$\Omega^{[2,0]}$, $\Omega^{[0,2]}$, and $\Omega^{[1,1]}$

In what follows we shall use only symmetric tensors. Thus the product
$dz\,d\bz$ will be understood as the symmetric one. This corresponds to the
embedding of $\Omega^{[1,1]}$ into  $T^*S\otimes T^*S\otimes\cc$ with symmetric image.

\item\label{pr3}  The complex conjugation in the bundles $(TS)^\cc$ and $(T^*S)^\cc$
 induces conjugacy bundle homomorphisms $\sigma: T^{(1,0)}S\to T^{(0,1)}S$ and
 $\sigma:\Omega^{(1,0)}S\to \Omega^{(0,1)}S$ with $\sigma^2=\id$. 

This $\sigma$ extends to conjugacy
 homomorphisms $\sigma: T^{[p,q]}S\to T^{[q,p]}S$ and $\sigma:\Omega^{[p,q]}S\to \Omega^{[q,p]}S$ for
 any pair of integers $p,q$. In the case $T^{[0,0]}=\Omega^{[0,0]}=\cc$ conjugacy
 $\sigma$ coincides with usual complex conjugation. We use the notation $\bar \alpha$
 for the complex conjugate of any section $\alpha$.

 If $z=x+\cplxi y$ is a local complex coordinate on $S$ as above, then
 $dz^p\,d\bz^q$ is a local trivialisation of $\Omega^{[p,q]}S$, and its complex
 conjugate is $dz^q\,d\bz^p$. 

\item In the special case $p=q=1$ we obtain that
 $dz\,d\bz=dx^2+dy^2$ is real.  Further, the local expression of the metric
 in the local coordinates $x,y$ is $g=\lambda(dx^2+dy^2)$. So we can consider $g$ as
 a global section of $\Omega^{[1,1]}S$. 

 The metric volume form $\omega=\lambda dx\land dy$ is the real section of the bundle
 $\Lambda^2S\otimes\cc$ of complex $2$-forms on $S$. Moreover, under natural isomorphisms
 $\Omega^{[1,1]}S\xrar{\,\cong\,}\Lambda^2S\otimes\cc$ the metric form $g=\lambda(dx^2+dy^2)$ is mapped
 onto the volume form $\omega=\lambda dx\land dy$. The difference between $\Omega^{[1,1]}S$ and
 $\Lambda^2S\otimes\cc$ lies in that how these bundles are embedded in $T^*S\otimes T^*S\otimes\cc$.

\item\label{pr5} Let $\nabla$ be the metric covariant derivative on $S$. Then it can be
 applied to any section $\alpha$ of any bundle $\Omega^{[p,q]}S$ (or $T^{[p,q]}S$) and
 the result lies in the product $\Omega^{[p,q]}S\otimes(T^*S)^\cc$ (or respectively
 $T^{[p,q]}S\otimes(T^*S)^\cc$) which is the sum $\Omega^{[p+1,q]}S\oplus\Omega^{[p,q+1]}S$
 (respectively $T^{[p-1,q]}S\oplus T^{[p,q-1]}$). We denote by $\nabla^{(1,0)}$ and
 $\nabla^{(0,1)}$ the component of $\nabla\alpha$ lying in $\Omega^{[p+1,q]}S$ and $\Omega^{[p,q+1]}S$,
 respectively. In ``dual'' notation, $\nabla^{(1,0)}$ and $\nabla^{(0,1)}$ act from
 $T^{[p,q]}S$ to $T^{[p-1,q]}S$ and $T^{[p,q-1]}$, respectively.

 Explicit calculation show that in coordinates $z=x+\cplxi y$ with
 $g=\lambda(dx^2+dy^2)$ one has
\[\textstyle\begin{split}
\nabla^{(1,0)}(fdz^pd\bz^q)= & \textstyle 
\big(\frac{\partial}{\partial z}f -p\frac{\partial\lambda}{\lambda\partial z}f\big)dz^{p+1}d\bz^q
\\
\nabla^{(0,1)}(fdz^pd\bz^q)= & \textstyle
\big(\frac{\partial}{\partial\bz}f -q\frac{\partial\lambda}{\lambda\partial\bz}f\big)dz^pd\bz^{q+1}
\end{split}
\]

In particular, $\nabla^{(1,0)}g= \nabla^{(1,0)}g=0$, $\nabla^{(0,1)}$ coincides with
$\frac{\partial}{\partial\bz}$ for sections of holomorphic bundles $T^{[p,0]}S$,
$\Omega^{[p,0]}S$, so that such a section $\alpha$ is holomorphic iff $\nabla^{(1,0)}\alpha=0$.

\item\label{curv-R}
 The operators $\nabla^{(1,0)}$ and $\nabla^{(0,1)}$ of complex and anti-complex
 covariant differentiation do not compute. We have  
\begin{equation}\label{commutator}
(\nabla^{(1,0)}\nabla^{(0,1)} - \nabla^{(0,1)}\nabla^{(1,0)})fdz^pd\bz^q = 
\frac{q-p}{2}\cdot R\cdot g\cdot fdz^pd\bz^q
\end{equation}
where $R$ is the Gauss curvature of the metric $g$. The formula is deduced as
follows. First, one computes the formula directly for the case $f\equiv1$. Then 
using Leibniz rule one shows that $(\nabla^{(1,0)}\nabla^{(0,1)} -
\nabla^{(0,1)}\nabla^{(1,0)})$ %
is not a differential operator but a bundle homomorphism, and hence the
formula holds for arbitrary $f$.

\item To simplify notation, we denote the complex and anti-complex covariant
 differentiations $\nabla^{(1,0)}$ and $\nabla^{(0,1)}$ by $(\cdot)_{;z}$ and
 $(\cdot)_{;\bz}$, respectively. Iterated differentiations are denoted like
 $(\cdot)_{;z\bz\bz z}$ and so on, the order of symbols $z$/$\bz$ coincides
 with the order of corresponding differentiations.
In particular, the formula \eqref{commutator} now reads
\begin{equation}\label{comm-z-bz}\begin{split}
(fdz^pd\bz^q)_{;z\bz} - (fdz^pd\bz^q)_{;\bz z} = 
\frac{p-q}{2}\cdot R\cdot g\cdot fdz^pd\bz^q;\\
(f\d_{z\mathstrut}^p\d_{\bar z\mathstrut}^q)_{;z \bz\mathstrut}\vph - 
(f\d_{z\mathstrut}^p\d_{\bar z\mathstrut}^q)_{;\bz z\mathstrut}\vph = 
\frac{q-p}{2}\cdot R\cdot g\cdot f\d_{z\mathstrut}^p\d_{\bar z\mathstrut}^q.
\end{split}
\end{equation}

\item \label{PoiBra-z} We shall use the formula
\begin{equation}\label{PoiBra-z-form}
\textstyle
f_{;z}h_{;\bz}-f_{;\bz}h_{;z} = 
\frac{\cplxi}{2}g\{f, h\}_g.
\end{equation}
for the Poisson bracket of functions $f$ and $h$ on $S$ with respect to the metric
symplectic form $\omega_g=\lambda dx\land dy=\frac{\cplxi}{2} \lambda dz\land d\bz$. Its calculation is as
follows:
\[
\textstyle
f_{;z}h_{;\bz}-f_{;\bz}h_{;z} = 
\frac{ df\land dh }{ dz\land d\bz } \,dz\,d\bz= 
\frac{\cplxi }{2}\frac{ df\land dh }{ \frac\cplxi2 {\cdot}\lambda dz\land d\bz } \lambda dzd\bz= 
\frac{\cplxi}{2}g\{f, h\}_g.
\]

\item Sections of bundles $\Omega^{[k,0]}S$, especially holomorphic ones, are
 called \slsf{$k$-differentials}. For this reason we call sections of bundles
 $T^{[k,0]}S$ {\it\;$k$-codifferentials}.
\end{enumerate}

\bigskip
Let us now recalculate the Poisson bracket $\{F,H\}$. Recall that we consider
the functions $F$ on $T^*S$ which are polynomial in momenta. In local
coordinates $z,\bz,p,\barp$, such a function is a polynomial in $p,\barp$
whose coefficients are smooth functions on $S$. With every such polynomial
$\sum_{ij=0}^df_{ij}p^i\barp^j$ we associate a section of the bundle sum
$\bigoplus_{i,j=0}^dT^{[i,j]}S$ using the rule $p\mapsto\frac{\partial}{\partial z},\barp\mapsto\frac{\partial}{\partial\bz}$
and extending it in the obvious way on polynomials. Let us denote this map by
$\bft:F\mapsto\bft(F)$. In particular, $\bft(H)=\frac12g\inv$.

\begin{lem}\label{braFH-1} $\bft(\{F,H\})=\bft(H)\cdot \nabla(\bft(F))$.
\end{lem}

\proof Let us observe that since both sides satisfy the Leibniz rule in
$F$ and commute with complex conjugation, it is sufficient to check the
formula in the case $F=p$. The rest follows. \qed

\medskip
In view of the lemma, we can replace the functions on $T^*S$ polynomial in
momenta by their $\bft$-images which are finite sums of sections of
$T^{[i,j]}S$.  In our special case first integrals cubic in momenta correspond
to sums $\half(A+\bar{A}+ B +\bar{B})=\Re(A+B)$
where $A$ is a section of $T^{[3,0]}S$ and $B$ is a section of
$T^{[2,1]}S$.  The equation $\Re(\nabla A+\nabla B)=0$ takes values in the sum
$\oplus_{j=-1}^3 T^{[2-j,j]}S$, and homogeneous components of the equation are:
$\nabla^{(0,1)}A=0$, $\nabla^{(1,0)}A+\nabla^{(0,1)}B=0$, the complex
conjugates of these two equations, and $\Re(\nabla^{(1,0)}B)=0$.  As we have seen
above, the equation $\nabla^{(0,1)}A=0$ on a section of $T^{[3,0]}S$ means its
holomorphicity.

\subsection{Component $B$ and the principle equation.} \label{B}
Let us now consider the last equation $\Re(\nabla^{(1,0)}B)=0$.  The tensor $B$ is a
section of $T^{[2,1]}S$, so $g^2\cdot B$ is section of the bundle
$\Omega^{[0,1]}S$. Hence $g^2\cdot B$ has the form $g^2\cdot B=\beta(x,y)\cdot d\bz$ for some
complex function $\beta(x,y)=\beta_1(x,y)+\cplxi{\cdot}\beta_2(x,y)$ with real and imaginary
components $\beta_1,\beta_2$.  By property (\ref{pr5}) on page \pageref{pr5},
\[\textstyle
g^2 \nabla^{(1,0)}B= \nabla^{(1,0)}(g^2B)= \frac{\partial\beta}{\partial z}\,dzd\bz.
\]
Since $dzd\bz$ is real, the equation $\Re(\nabla^{(1,0)}B)=0$ is equivalent to
\[\textstyle
\Re(\frac{\partial\beta}{\partial z}) = \frac{\partial\beta_1}{\partial x} + \frac{\partial\beta_2}{\partial y}=0. 
\]
This equation can be seen as the closedness of the $1$-form $-\beta_2dx+\beta_1dy$.
Hence the equation $\Re(\nabla^{(1,0)}B)=0$ is equivalent to a local existence
of some \emph{real} function $K(x,y)$ such that $\beta_1=\frac{\partial K}{\partial y}$ and
$\beta_2=-\frac{\partial K}{\partial x}$. The latter two conditions can be rewritten as a single
complex equation $\beta= -2\cplxi \frac{\partial K}{\partial\bz}$, or equivalently
\begin{equation}\label{B=barzK}
B=-2\cplxi g^{-2}K_{;\bz}.
\end{equation}
The latter equation is identical with the equations (2.7), (2.8) in
\cite{Du-Ma-To}.

\smallskip%
Substituting \eqref{B=barzK} in the remaining equation
$\nabla^{(1,0)}A+\nabla^{(0,1)}B=0$, we obtain the following necessary and
 sufficient condition:

\begin{lem}\label{princ-eq} Given a metric $g$ and a
 holomorphic $3$-codifferential $A$ on a surface $S$, $A$ is compatible with $g$ a
if and only if the equation
\begin{equation}\label{Eqz}
\textstyle
\cale_z:\ \   K_{;\bz\bz}= -\frac{\cplxi}{2} g^2 A_{;z}
\end{equation}
has a smooth \emph{real-valued} solution $K$,  and in this case $F=\Re(A+B)$ with
 $B=-2\cplxi g^{-2} K_{;\bz}$ is a cubic integral with the
 Birkhoff-Kolokoltsov tensor   $A$. 
\end{lem} 

 We call \eqref{Eqz} the \slsf{principle equation.}

\section{Proof of Theorem 1.1.}

\subsection{Calculation of further equations.}
We apply the prolongation-projection method to the equation \eqref{Eqz}
considering the function $K$ as an unknown function and $g$, $A$ (and their
derivatives) as known parameters.  The equations obtained by
prolongation-projection procedure are called \slsf{deduced}.  Instead of
considering the real and imaginary parts of equations, we shall mostly use
complex equations and their complex conjugates. In particular, the complex
conjugate to \eqref{Eqz} is
\begin{equation}\label{Eqbz}
\textstyle
K_{;zz}=   \frac{\cplxi}{2} g^2  \bar{A}_{;\bz}
\end{equation}
Differentiating \eqref{Eqz} and \eqref{Eqbz} with respect to $z$ and $\bz$
(first prolongation), we obtain the following 4 equations:
\begin{equation}\label{EqKzzz}
\textstyle
K_{;zzz}=  \frac{\cplxi}{2} g^2 \bar{A}_{;\bz z} \;\;\;\;
K_{;zz\bz}=  \frac{\cplxi}{2} g^2 \bar{A}_{;\bz\bz}\;\;\;\;
K_{;\bz\bz z}=  -\frac{\cplxi}{2} g^2 A_{;zz}\;\;\;\;
K_{;\bz\bz\bz}=  -\frac{\cplxi}{2} g^2 A_{;z\bz}
\end{equation}
This gives us the expression of all partial derivatives of $K$ of order $3$
via lower order derivatives, which are ``hidden'' due to the covariant form of
the equations and do not appear explicitly.  In particular, the system
 $(\cale_z,\cale_\bz)$ has finite type.

 Differentiating once more (second prolongation), we obtain 6 equations of
 order $4$: 3 as order-2 derivatives of \eqref{Eqz} and 3 more as derivatives
 of its complex conjugates:
\begin{equation}\label{EqKz4}
\textstyle
\begin{split}\textstyle
K_{;zzzz}= \ \ \,\frac{\cplxi}{2} g^2 \bar{A}_{;\bz zz} \qquad
K_{;zzz\bz}= \ \ \,\frac{\cplxi}{2} g^2  \bar{A}_{;\bz z\bz} \qquad
K_{;zz\bz\bz}= \ \ \,\frac{\cplxi}{2} g^2  \bar{A}_{;\bz \bz \bz} 
\\\textstyle
K_{;\bz\bz zz}= -\frac{\cplxi}{2} g^2  A_{;zzz}\qquad
K_{;\bz\bz z\bz}= -\frac{\cplxi}{2} g^2  A_{;zz\bz}\qquad
K_{;\bz\bz\bz\bz}= -\frac{\cplxi}{2} g^2  A_{;z\bz\bz}
\end{split}
\end{equation}
Since there are 5 partial derivatives of $K$ of order $4$, at this
prolongation-projection  step  we can obtain one new equation.

Comparing $\Omega^{[p,q]}$-types of expressions in \eqref{EqKz4} we see that the
desired equation should be the difference of equations %
$K_{;zz\bz\bz}=  \frac{\cplxi}{2} g^2 \bar{A}_{;\bz \bz \bz}$ and %
$K_{;\bz\bz zz}= -\frac{\cplxi}{2} g^2 A_{;zzz}$:
\begin{equation}\label{Ersh1}
\textstyle
(K_{;zz\bz\bz} - K_{;\bz\bz zz}) - 
 \frac{\cplxi}{2} g^2(\bar{A}_{;\bz \bz \bz} + A_{;zzz}) 
\end{equation} 
Computing the difference $K_{;zz\bz\bz} - K_{;\bz\bz zz}$ we use the formula
for commutator of covariant derivatives, see property \eqref{curv-R}
on page \pageref{curv-R}. This gives
\[
\begin{split}
&K_{;zz\bz\bz} - K_{;\bz\bz zz} = 
\\&
(K_{;zz\bz\bz} - K_{;z\bz z\bz})_1 + 
(K_{;z\bz z\bz} -K_{;z\bz \bz z})_2 +
(K_{;z\bz \bz z}- K_{;\bz z \bz z})_3+
(K_{;\bz z \bz z}- K_{;\bz\bz zz})_4
=
\\&\textstyle
((\half R g K_{;z})_{;\bz})_1 - 
((\half R g K_{;\bz} )_{;z})_4= 
\frac{g}{2} (K_{;z} R_{;\bz}-K_{;\bz}R_{;z}).
\end{split}
\]
Here we use \eqref{commutator} which shows that the differences $(\ldots)_2$ and
$(\ldots)_3$ vanish and gives the formulas above for the differences $(\ldots)_1$ and
$(\ldots)_4$.  So using \eqref{PoiBra-z-form} we can write down the first deduced
equation (new equation in the terminology of §\ref{prol1}).
\begin{equation}\label{Ersh}
 \{K,R\}= 4\Re(A_{;zzz})=:D_0
\end{equation}
where we denote the right hand side by
$D_0=:4\Re(A_{;zzz})= 2A_{;zzz}+2\bar{A}_{;\bz\bz\bz}$.

\medskip \state Remark.  We see that the first deduced equation does not
contain second derivatives of $K$. This is not a coincidence but follows from
the condition that the higher coefficient in the \eqref{Eqz} is constant, see
\cite{Kr-Ly} for details.

At this point we can solve the problem of existence of cubic integrals 
with a given Birkhoff-Kolokoltsov tensor $A$ for metrics with constant
curvature. 

\begin{prop}\label{R=const} Let $g$ be a metric metric on a surface $S$ of
 constant curvature and $A$ a holomorphic $3$-codifferential. Then $A$ is
 compatible with $g$ if and only if $A$ satisfy the PDE $D_0=0$.

 Moreover, for a given point $P\in S$ there exists a unique cubic integral of
 the form $F=\Re(A)+2H\cdot b^ip_i$ with a vector field $\vec b:=b^i\partial_{x^i}$ having
 prescribed values $\vec b(P)=\vec v\in T_PS$ and
 $\mathsf{rot}\,\vec{b}(P)=w\in\rr$ at the point $P$. 
\end{prop}

\proof Let $(x^1,x^2)$ be any given coordinate system on $S$. The
formula \eqref{EqzTensor} gives the covariant tensor form of the principle
equation \eqref{Eqz}. Denote by $\cale_{ij}$ the components of this tensor
form.  Let $\cals$ be the system consisting of equations $\cale_{ij}$ and
their 1st order derivatives $\cale_{ij;k}$.  Introduce new variables
$\bfk_i=\bfk_i(x^1,x^2),\bfk_{ij}=\bfk_{ij}(x^1,x^2)$, substitute the latter
into equations $\cale_{ij}$ instead of the corresponding derivatives
$K_{;ij}$, and denote the obtained inhomogeneous linear-algebraic equations by
$\cale_{ij}(\bfk)$. Consider the system $\cals(\bfk)$ consisting of equations
$\cale_{ij}(\bfk)$, their partial derivatives $\cale_{ij}(\bfk)_{;k}$, the
equations $(\bfk_{i})_{;j}=\bfk_{ij}$ and the consistency equations
$\bfk_{12}=\bfk_{21}$ and $(\bfk_{ij})_{;k}=(\bfk_{ik})_{;j}$. Then the
systems $\cals$ and $\cals(\bfk)$ are equivalent: the inverse substitution of
$K_{;i},K_{;ij}$ instead of $\bfk_i,\bfk_{ij}$ transforms $\cals(\bfk)$ into
$\cals$ (turning consistency equations into identities). Moreover, both system
are involutive simultaneously.

Now, since the system $\cals(\bfk)$ consists of 1st order linear PDEs and is
involutive, it admits a local solution by the classical Frobenius theorem (see
eg.\ \cite{Iv-La}, §\,1.9). Moreover, for every $P\in S$ and every values
$\bfk_i(P),\bfk_{ij}(P)$ satisfying \emph{algebraic} equations
$\cale_{ij}(\bfk)$ at the point $P$ there exists a local solution
$(\bfk_i(x),\bfk_{ij}(x))$. For this $\bfk$-solution there exists a function
$K(x)$, unique up to adding a constant, which satisfies the relations
$K_{;i}=\bfk_i,K_{;ij}=\bfk_{ij}$. 

To calculate the number of parameters, we use linear-algebraic equations
$\cale_{ij}(\bfk)$ and express the functions $\bfk_{12}=\bfk_{21}$ and
$\bfk_{22}$ via $\bfk_{1},\bfk_{2}$, and $\bfk_{11}$. Substituting these
formulas in the remaining PDEs, we reduce the number of unknown functions to
three. Let $\cals'(\bfk)$ be the obtained system.  Notice that some of the substitutions
above mean the application of prolongation-projection procedure, for
example, substitution of the formulas $\partial_i\bfk_{jk}=\calf_{ijk}$ in equation
$(\bfk_{ij})_{;k}=(\bfk_{ik})_{;j}$.  However, the involutivity of
$\cals(\bfk)$ implies that every equation in $\cals'(\bfk)$ is still a
\emph{linear-algebraic} consequence of the system $\cals(\bfk)$. In turn, this
implies the involutivity of the system $\cals'(\bfk)$: Indeed, every
differential consequence of $\cals'(\bfk)$ is a differential consequence of
$\cals(\bfk)$, and hence is a linear-algebraic consequence.

By the construction, the system $\cals(\bfk)$ is linear-algebraic equivalent
to a system containing the equations of the form
$\partial_i\bfk_j=\calf_{ij},\partial_i\bfk_{jk}=\calf_{ijk}$ where $\calf_{ij},\calf_{ijk}$
are linear inhomogeneous expressions in $\bfk_i,\bfk_{ij}$. It follows that
the system $\cals'(\bfk)$ is also linear-algebraic equivalent to a system
\emph{consisting of} the equations of the form
$\partial_i\bfk_j=\calf'_{ij},\partial_i\bfk_{11}=\calf'_{i11}$ where
$\calf'_{ij},\calf'_{i11}$ are as above linear inhomogeneous in
$\bfk_i,\bfk_{11}$, and we may assume that the system $\cals'(\bfk)$ is of
this form. 

Now the generic solution of the system $\cals'(\bfk)$ is constructed as
follows: Given a point $P$ and values $\bfk_{1}(P),\bfk_{12}(P),\bfk_{11}(P)$
at $P$, we fix local coordinates $x:=x^1$ and $y:=x^2$,  integrate
the ODE system $(\partial_{x}\bfk_i=\calf_{1i},\partial_{x}\bfk_{11}=\calf_{111}$ along some
interval on the $x$-axis, and then integrate
the ODE systems $(\partial_{y}\bfk_i=\calf_{2i},\partial_{y}\bfk_{11}=\calf_{211}$ along 
every interval parallel to $y$-axis.

It remains to show that the prescribed values as in the assertion of the
proposition can be used to parametrised a general solution.  Indeed, as such
parameters we can use the values of the 1st order derivatives of $K$ and the
Laplacian $\Delta_gK$ (instead of $K_{;11}$ at a given point $P\in S$. Now recall
that by §\,\ref{B}, $b^i$ and $K$ are related as $K_{;i}=\omega_{ij}b^j$ so that
$\vec b$ is the Hamiltonian vector field on $S$ of the (Hamiltonian) function
$K$ (with respect to the form $\omega_g$). Consequently, $\vec b$ is the
skew-gradient of $K$ and the Laplacian $\Delta_gK$ is the rotor
$\mathsf{rot}_{\,}\vec{b}$.
\qed

\subsection{Proof of Theorem 1.1.}
\label{proof-th1} %
Let us now turn to the (most interesting) general case when $R$ is
non-constant.  Write the equation \eqref{Ersh} as
$K_{;z}R_{;\bz}-K_{;\bz}R_{;z}=\frac{\cplxi}{2} g D_0$.  Differentiating it in
$z$ and $\bz$ we obtain two more equations of order 2 in $K$:
\begin{equation}\label{Ershz}
\textstyle
\begin{split} \textstyle
K_{;zz}R_{;\bz} +K_{;z}R_{;\bz z} 
-K_{;\bz z}R_{;z}-K_{;\bz}R_{;zz}= \frac{\cplxi}{2} g{\cdot}(D_0)_z
\\ \textstyle
K_{;z\bz}R_{;\bz} +K_{;z}R_{;\bz \bz} 
-K_{;\bz \bz}R_{;z}-K_{;\bz}R_{;z\bz}= \frac{\cplxi}{2} g{\cdot}(D_0)_\bz
\end{split} 
\end{equation}
Together with \eqref{Eqz} and its conjugate we now have 4 equations of second
order on $K$. So we expect to obtain at this step one new 
equation of order 1 on $K$. It is obtained as follows:

Add equations \eqref{Ershz} with coefficients $R_{;\bz}$ and $R_{;z}$:
\[
\begin{split} 
&K_{;zz}R_{;\bz}R_{;\bz}  
     +K_{;z}R_{;\bz z}R_{;\bz}
-K_{;\bz}R_{;zz}R_{;\bz}+
\\ 
&
\textstyle
K_{;z}R_{;\bz \bz}R_{;z} 
-K_{;\bz \bz}R_{;z}R_{;z}
     - K_{;\bz}R_{;z\bz}R_{;z}
= \frac{\cplxi}{2} g{\cdot}(D_0)_{;z}  R_{;\bz}  
+ \frac{\cplxi}{2} g{\cdot}(D_0)_{;\bz} R_{;z}
\end{split} 
\]
Rearranging it, we obtain
\[
\begin{split} 
&
K_{;z}(R_{;z}R_{;\bz})_{;\bz} 
-K_{;\bz}(R_{;z}R_{;\bz})_{;z}
\textstyle
= \frac{\cplxi}{2} g{\cdot}(D_0)_{;z}R_{;\bz}  
     + \frac{\cplxi}{2} g{\cdot}(D_0)_{;\bz} R_{;z}
+K_{;\bz \bz}R_{;z}R_{;z}-K_{;zz}R_{;\bz}R_{;\bz} .
\end{split} 
\]
The expression $R_{;z}R_{;\bz}$ equals $\frac{1}{4}g|\nabla R|^2$. Let us denote
$\frac12|\nabla R|^2$ by $\varphi_1$. So the left hand side of the latter equation is
$\frac{\cplxi}{4} g^2\{K,\varphi_1\}$.  On the right hand side we have 
 $\frac{\cplxi}{2} g\big( (D_0)_{;z}R_{;\bz}+(D_0)_{;\bz}R_{;z}\big)$ which
 equals $\frac{\cplxi}{4} g^2\langle\nabla D_0,\nabla R\rangle$. In the last two terms we use
\eqref{Eqz}. So finally we obtain the equation
\begin{equation}\label{Ezh}
 \{K,\varphi_1\}=  \langle\nabla D_0,\nabla R\rangle -4\Re(A_{;z}\cdot(R_{;z})^2)=:D_1,
\end{equation}
in which we denoted the right hand side by $D_1$.

\medskip
From this points we can produce many equations, each of the form $\{K,\varphi\}=D$
where $\varphi$ and $D$ are real functions which are certain  differential
expressions involving the metric $g$ and the tensor $A_{;z}$. In fact, we have
two such ``dummy'' procedures:

The first one is the formal repetition of the deduction
\eqref{Ersh}$\Longrightarrow$\eqref{Ezh} and uses the fact that the explicit expressions for
$R$ and $D_0$ were not involved. So from any equation $\{K,\varphi\}=D$ we can obtain
a new equation $\{K,\varphi'\}=D'$ with $\varphi':=\half|\nabla\varphi|^2$ and $D':=\langle\nabla
 D,\nabla\varphi\rangle-4\Re(A_{;z}\varphi_{;z}^2)$.

The second procedure is based on the Jacobi identity for Poisson brackets:
starting with any two equations of the form $\{K,\varphi'\}=D'$ and $\{K,\varphi''\}=D''$, we
obtain
\begin{equation}\label{dummy}
\{K,\{\varphi',\varphi''\}\}=\{\{K,\varphi'\},\varphi''\} - \{\{K,\varphi''\},\varphi'\}= \{D',\varphi''\}-\{D'',\varphi'\} 
\end{equation}
and we can set $\varphi''':= \{\varphi',\varphi''\}$ and $D''':=\{D',\varphi''\}-\{D'',\varphi'\}$.

Let us notice that both procedures are application of prolongation-projection
Method: In the first case we differentiate the equation $\{K,\varphi\}=D$ in $z$ and
$\bz$, add the equation \eqref{Eqz} and its complex conjugate, and make
linear-algebraic manipulations on the system of 4 equations excluding second
order derivatives of $K$.  In the second case we differentiate both equations
$\{K,\varphi'\}=D'$ and $\{K,\varphi''\}=D''$ and then make the same linear-algebraic
manipulations on the system of 4 equations.

To unify the notation, we set
\begin{equation}\label{phi-D-0-3}
\textstyle
\begin{split}
&\textstyle
\varphi_0:=R,\qquad
\varphi_1:=\half|\nabla\varphi_0|^2,\qquad
\varphi_2:=\{\varphi_0,\varphi_1\},\qquad
\varphi_3:=\frac12|\nabla\varphi_1|^2,
\\
&
\textstyle
D_0:=4\,\Re(A_{;zzz}),\qquad\qquad\qquad\quad\!
D_1:=\langle\nabla D_0,\nabla\varphi_0\rangle-4\,\Re\big(A_{;z}\cdot((\varphi_0)_{;z})^2\big),
\\
&
\textstyle
D_2:=\{D_0,\varphi_1\}-\{D_1,\varphi_0\},
\qquad\;
D_3:=\langle\nabla D_1,\nabla\varphi_1\rangle-4\,\Re\big(A_{;z}\cdot((\varphi_1)_{;z})^2\big).
\end{split}
\end{equation}%
Thus the equation $\{K,\varphi_2\}=D_2$ is obtained from $\{K,\varphi_0\}=D_0$ and
$\{K,\varphi_1\}=D_1$ using Jacobi identity, and the equation $\{K,\varphi_3\}=D_3$ from
$\{K,\varphi_1\}=D_1$ using \eqref{Ersh} and its complex conjugate. 

\smallskip%
It appears that the system consisting of 4 equations $\cale_i$:
$\{K,\varphi_i\}=D_i$, $i=0,\ldots,3$ is involutive and differentially equivalent
 to the original equation $\cale_z$. We state a more general property which
 will be used also in the proof of Theorem \ref{th:main2}.

\begin{lem}\label{lin-alg} 
 \smallskip\noindent%
 {\bf(a)} Let the coefficients of the equations $\cale'=\{K,\varphi'\}-D'$,
 $\cale''=\{K,\varphi''\}-D''$, and $\cale'''=\{K,\varphi'''\}-D'''$ are related as $\varphi''':=
 \{\varphi',\varphi''\}$ and $D'':=\{D',\varphi''\}-\{D'',\varphi'\}$. Then the equation $\cale'''$ is a
 \emph{linear-algebraic} consequence of the equations
 $\cale'_{;x},\cale'_{;y}$,
 $\cale''_{;x},\cale''_{;y}$.

\smallskip\noindent%
{\bf(b)} Let the coefficients of the equations
 $\cale'=\{K,\varphi\}-D$ and $\cale^*=\{K,\varphi^*\}-D^*$ are related as $\varphi^*:=\half|\nabla\varphi'|^2$
 and $D^*:=\langle\nabla D',\nabla\varphi'\rangle-4\,\Re(A_{;z}(\varphi'_{;z})^2)$. Then the equation $\cale^*$ is a
 \emph{linear-algebraic} consequence of the equations
 $\cale'_{;x},\cale'_{;y}$, and the equations
 \begin{enumerate}
 \item[(I)] $\cale_z$ and its conjugate $\cale_\bz$. %
  \suspend{enumerate}

\smallskip\noindent%
{\bf(c)} Under hypotheses of parts {\bf(a)} and {\bf(b)} assume
additionally that $\{\varphi',\varphi''\}=\varphi'''$ is non-vanishing at a generic point. Let
also the equation $\cale^\dag=\{K,\varphi^\dag\}-D^\dag$ be defined by $\varphi^\dag:=\half|\nabla\varphi''|^2$
and $D^\dag:=\langle\nabla D'',\nabla\varphi''\rangle-4\,\Re(A_{;z}(\varphi''_{;z})^2)$.  Then the equations {\rm(I)}
are \emph{linear-algebraic} consequences of the following set of equations:
\resume{enumerate}
\item[(II)] equations $\cale''',\cale^*,\cale^\dag$;
\item[(III)] the first order derivatives $\cale'_{;x}$,
 $\cale'_{;y}$, $\cale''_{;x}$, $\cale''_{;y}$.
\end{enumerate}
\end{lem}

Recall that ``linear-algebraic consequence'' means that new equations appear
as linear combinations of old ones with coefficients which are rational
functions of coefficients of the old equations.

\proof Clearly, parts {\bf(a)} and {\bf(b)} are simply a
restating of the ``dummy'' procedures introduced above.

\smallskip\noindent%
{\bf(c)} Rewrite the equation $\cale'$ and $\cale''$ in the form
$K_z\varphi'_{;\bz}-K_\bz\varphi'_{;z}=2\cplxi g D'$ and
$K_z\varphi''_{;\bz}-K_\bz\varphi''_{;z}=2\cplxi g D''$. %
Then the determinant of the linear system is
$\varphi'_{;\bz}\varphi''_{;z}-\varphi''_{;\bz}\varphi'_{;z}=-2\cplxi g\cdot\{\varphi',\varphi''\}$, and the solution is
given by the formulas $K_i =\calf_i$ with $\calf_i:=\frac{1}{\{\varphi',\varphi''\}}\cdot \det
\smpmatr{D'\ & -\varphi'_{;x^i}\\ D''\ &-\varphi''_{;x^i}}$. %
Differentiating the formulas, we conclude that there exists a \emph{unique}
solutions of the equations (III), considered as a system of linear equations
on $K_{;x^ix^j}$, and this solution is given by
${K_{;x^ix^j}}=(\calf_i)_{;x^j}$. In particular,
$(\calf_1\vph)_{;x^2}=(\calf_2\vph)_{;x^1}$. Notice that each
$(\calf_i)_{;x^j}$ has the form $\frac{\wt\calf_{ij}}{\{\varphi',\varphi''\}^2}$ for some
PDOs $\wt\calf_{ij}$ which are \emph{polynomial} in $\varphi',\varphi'',D',D''$ and their
derivatives. This shows the fact that the determinant of the matrix of leading
coefficients of the system (III) equals $\{\varphi',\varphi''\}^2$ (up to sign). Indeed, the
solutions $K_{;x^ix^j}=(\calf_i)_{;x^j}$ can be obtained by means of
linear-algebraic manipulations with equations (III).

The key point in the proof is that we can obtain the equations (II)
substituting the expressions $K_{;x^ix^j}=(\calf_i\vph)_{;x^j}$ in the
equations (I). This is also linear-algebraic
manipulations with equations. The assertion {\bf(c)} claims that this
operation is invertible, which is true.
\qed

\medskip%
Define the expressions $\calg_i$ using the following
Jacobi-like expressions:
\[
\begin{split}
\calg_0&= \{\varphi_1,\varphi_2\}\cdot D_3+\{\varphi_2,\varphi_3\}\cdot D_1+\{\varphi_3,\varphi_1\}\cdot D_2,
\\
\calg_1&= \{\varphi_0,\varphi_2\}\cdot D_3+\{\varphi_2,\varphi_3\}\cdot D_0+\{\varphi_3,\varphi_0\}\cdot D_2,
\\
\calg_2&= \{\varphi_0,\varphi_1\}\cdot D_3+\{\varphi_1,\varphi_3\}\cdot D_0+\{\varphi_3,\varphi_0\}\cdot D_1,
\\
\calg_3&= \{\varphi_0,\varphi_1\}\cdot D_2+\{\varphi_1,\varphi_2\}\cdot D_0+\{\varphi_2,\varphi_0\}\cdot D_1.
\end{split} 
\]
Up to normalisation, the  expressions $\calg_i$ are
the $(3×3)$-minors of the extended matrix of the coefficients of the equations
$\cale_0,\ldots,\cale_3$ with $i$-th row excluded. In particular,
\[\textstyle
\calg_2 := \frac{\cplxi}{2g}\det
\pmatr{(\varphi_0)_{;z} & (\varphi_0)_{;\barz} & D_0\\ 
(\varphi_1)_{;z} & (\varphi_1)_{;\barz} & D_1\\ 
(\varphi_3)_{;z} & (\varphi_3)_{;\barz} & D_3
}
\qquad\text{and}
\qquad
\calg_3 := \frac{\cplxi}{2g}\det
\pmatr{(\varphi_0)_{;z} & (\varphi_0)_{;\barz} & D_0\\ 
(\varphi_1)_{;z} & (\varphi_1)_{;\barz} & D_1\\ 
(\varphi_2)_{;z} & (\varphi_2)_{;\barz} & D_2
}
.
\]
Every expression $\calg_i$ is a PDE on the metric $g$ and the
$3$-codifferential $A$.

\begin{corol}\label{complete}  Assume that $\varphi_2$ is n on-vanishing. 
 Then the system $\cals:=(\cale_0,\cale_1)$ is involutive if and only if the
 integrability condition $\calg_3=0$ is fulfilled. A solution $K$ of this
 system satisfies the equation $K_{;\bz\bz}=-\frac{\cplxi}{2}g^2A_{;z}$ if and
 only if the integrability condition $\calg_2=0$ is fulfilled.
\end{corol}

At this point we can give the proof of Theorem \ref{th:main1}.

\begin{prop}\label{prop-thm1} Let $g=\lambda(x,y)(dx^2+dy^2)$ be a metric and $A$ a
 holomorphic $3$-codifferential such that the differential expression
 $\varphi_2=\{\varphi_0,\varphi_1\}$ is non-vanishing.

 Then $g$ admits a cubic integral of the form $F=\Re(A+b\,p^2\bar p)$
 if and only if $\lambda$ satisfies the covariant PDEs
 $\calg_3,\calg_2$. Moreover, in this case  the
 complex valued function $b=b_1+\cplxi b_2$ is given by the formulas
 $b_1=\lambda^{-2}\calk_2, b_2=-\lambda^{-2}\calk_1$ with
\begin{equation}\label{K-i}
\calk_i:=\frac{1}{\{\varphi_0,\varphi_1\}}\cdot \det
\pmatr{D_0\  & -(\varphi_0)_{;x^i} \\ D_1\ & -(\varphi_1)_{;x^i}}.
\end{equation}
\end{prop}

\proof  By Lemma \ref{princ-eq} from
 §\,\ref{B}, the existence of a cubic integral $F$ of the form above is
 equivalent to the existence of a solution $K$ of the principle equation
 $K_{;\bz\bz}=-\frac{\cplxi}{2}g^2A_{;z}$, and in this case $F=\Re(A+B)$ with
 $B:=b\,p^2\bar p=-2\cplxi g^{-2} K_{;\bz}$ is the desired cubic integral.

As we have shown above, equations $\cale_0,\cale_1$ are differential
consequences of the real and imaginary parts of  $\cale_z$. So by 
Corollary \ref{complete},  equations  $\calg_2=0,\calg_3=0$
are necessary and sufficient conditions for solvability 
of $\cale_z$. The formulas \eqref{K-i} give simply the solution of a linear
system $(\cale_0,\cale_1)$.
\qed

Let us observe that under condition of non-vanishing of all brackets
$\{\varphi_i,\varphi_j\}$ any two of the equations $\calg_0,\ldots,\calg_3$ are linear-algebraic
consequence of remaining two.

\section{Metrics with $\{R, |\nabla R|^2_g\}_g^{\,}\equiv0$ and proof of Theorem
 1.2.}\label{degen}

In this section we consider the degenerate case when the function $\varphi_2$
introduced in \eqref{phi2} vanishes identically. This will be a general
assumption in this section unless the opposite is stated explicitly.  In this
case we additionally suppose that the Gauss curvature $R=\varphi_0$ is
non-degenerate, ie., the gradient of $R$ is non-vanishing.

\subsection{Proof of Theorem 1.2.}\label{thm2-proof}
Recall that the condition $\varphi_2\equiv0$ means that the Gauss curvature $R$ and the
square of its gradient $|\nabla R|^2$ are functionally dependent,
$\half|\nabla R|^2=f(R)$ for some function $f(r)$ of one variable.  Notice that the
function $\varphi_3=\half|\nabla\varphi_1|^2$ is also functionally dependent, since
\[\textstyle
\half|\nabla\varphi_1|^2=\half|\nabla f(R)|^2=(f'(R))^2\,\half|\nabla R|^2=(f'(R))^2f(R).
\]

\smallskip%
As in the non-degenerate case, we want to apply the prolongation-projection
method. Our starting system of equations is again
 $\cals_0:=\{\cale_z,\cale_\bz\}$.  Let us analyse the ``forking point'' in the
procedure. Since we assume non-vanishing of $\nabla R$, there will be no divergence
in results until Corollary \ref{complete}.  At this point the
degenerate case $\varphi_2\equiv0$ differs from the non-degenerate one in two aspects.

The first one  is that we must replace $\cale_1$ by another
equation. We shall find such an equation in a moment. The second aspect is
that in the case $\varphi_2\equiv0$ the result of the projection procedure is slightly
different as before.  Namely, now the system $(\cale_0,\cale_1)$ is degenerate
at every point since the left hand sides of the equations, $\{K,\varphi_0\}$ and
$\{K,\varphi_1\}$, are Lie derivatives of $K$ along proportional vector fields. Since
the equation $\cale_1$ still must be satisfied, we obtain a new differential
condition on $g$ and $A$: the right hand sides of $\cale_0,\cale_1$ must be
proportional with the same coefficient as the left ones. Clearly, this
condition is simply the application of the projection procedure to the system
$(\cale_0,\cale_1)$.

The resulting equation reads:
\begin{equation}\label{D-eq}
\cald_z := (\varphi_0)_{;z}D_1-(\varphi_1)_{;z}D_0 =0.
\end{equation}
Notice that even if this is formally a complex-valued equation and so two real
ones, the real and imaginary parts are equivalent under the condition
$\varphi_2\equiv0$. This means that the system of conditions  $\{\varphi_2,\cald_z\}$ contains only two
independent conditions. 

In a more explicit form the equation $\cald_z$ reads:
\begin{equation*}\textstyle
\cald_z =
R_{;z}\cdot\langle\nabla D_0,\nabla R\rangle  - R_{;z}\cdot4\cdot\Re(A_{;z}\cdot(R_{;z})^2)
- D_0\cdot(\half|\nabla R|^2)_{;z}.
\end{equation*}
Another possible form for this  condition is the differential $1$-form 
\[\textstyle
\det
{\scriptstyle
\left(\begin{matrix}
 d\varphi_0 & D_0 \\ d\varphi_1 & D_1
\end{matrix}\right)}
= \cald_zdz +\cald_\bz d\bz
= \cald_xdx +\cald_y dy
,
\]
where $\cald_\bz,\cald_x,\cald_y$ are defined in the obvious way.  
One can use any of the equation $\cald_x,\cald_y$ instead of $\cald_z$ provided the
corresponding derivative $\frac{\partial\varphi_0}{\partial x^i}$ is non-vanishing.

\smallskip%
Now we seek for the substitute for the equation $\cale_1$.  Since our previous
step was the projection, the next one is the prolongation of the equations.
Let $\cals$ be the system of equations obtained so far. Those are
$\cale_z,\cale_\bz$ and their derivatives up to order $2$, the equation
$\cale_0$ and its derivative $(\cale_0)_{;x}$, and the equations $\varphi_2$ and
$\cald_z$.

As the next step we are going to add to $\cals$ the second order derivatives
$(\cale_0)_{;x^ix^j}$.  Here we make the following easy observation.  It
follows from Lemma \ref{lin-alg} and the above consideration that the
equations $(\cale_0)_{;y}$ is a linear-algebraic combination of the equations
from $\cals$ of order $\leq2$ in $K$ and the condition $\cald_z$.
Similarly, $(\cale_0)_{;x}$ can be obtained as a linear algebraic combination
of the condition $\cald_z$ and the equations from $\cals$ of order $\leq2$ in
$K$, in which $(\cale_0)_{;x}$ is replaced by $(\cale_0)_{;y}$. Thus replacing
$(\cale_0)_{;x}$ by $(\cale_0)_{;y}$ we obtain a system equivalent to $\cals$
provided $\cald_z\equiv0$.

Consequently, under condition $\cald_z\equiv0$ the derivative
$(\cale_0)_{;yy}$ is a linear algebraic combination of the same equations and
its derivatives in $y$. In particular, $(\cale_0)_{;yy}$ is a linear algebraic
combination of equations from $\cals$ and the equation
$(\cale_0)_{;xy}$. Interchanging $x$ and $y$ we conclude that
$(\cale_0)_{;xx}$ is also a linear algebraic combination of equations from
$\cals$ and the equation $(\cale_0)_{;xy}$. In fact, each of the equations
$(\cale_0)_{;x^ix^j}$ is linear-algebraic equivalent to each other modulo the
system $\cals$ and the condition $\cald_z$. Thus adding
to $\cals$ an arbitrary single derivative $(\cale_0)_{;x^ix^j}$ instead all
three we obtain an equivalent new system. As such a derivative we choose
$(\cale_0)_{;xy}=(\cale_0)_{;z\bz}$.

\smallskip%
Now we make explicit calculation. The derivation $(\cale_0)_{;z\bz}$ gives
\begin{equation}\label{Ersh-z-bz}
\begin{array}{rcc} 
K_{;zz\bz}R_{;\bz} + K_{;zz}R_{;\bz\bz} +
K_{;z\bz}R_{;\bz z} + K_{;z}R_{;\bz z\bz}  && 
\\
-K_{;\bz z\bz}R_{;z}-K_{;\bz z}R_{;z\bz}
-K_{;\bz \bz}R_{;zz}-K_{;\bz}R_{;zz\bz}
&= &2\cplxi g\cdot(D_0)_{;z\bz}.
\end{array}
\end{equation}
Simplifying the obtained equation we apply the following relations: the rule
(\ref{curv-R}), the homolorphicity equations $A_{;\bz}=0, \bar A_{;z}=0$, and
the substitutions \eqref{Eqz}, \eqref{EqKzzz}. Let us notice that the latter
case we exclude higher order derivatives of $K$ and so apply the projection
procedure. Calculating, we obtain $K_{;zz\bz} = 
\frac{\cplxi}{2}g^2\bar A_{;\bz\bz}$ and 
\[
\textstyle
K_{;\bz z\bz}=K_{;\bz\bz z}+ (K_{;\bz z\bz}-K_{;\bz\bz z})
=-\frac{\cplxi}{2}g^2 A_{;zz} -  \frac{1}{2}gRK_{;\bz}.
\]
Besides we also use the relations 
\[\textstyle
R_{;z\bz}=R_{;\bz z}= \frac{1}{4} g\Delta R
\qquad
\text{and}
\qquad
R_{;zz\bz} = R_{;z\bz z} + \frac{1}{2}gRR_{;z}
 = \frac{1}{4}g \Delta R_{;z} + \frac{1}{2}gRR_{;z}
\]
where $\Delta=\Delta_g$ is the metric Laplace operator.
Substitution of these relations yields
\begin{equation*}
\begin{split} \textstyle
\frac{\cplxi}{2}g^2\bar A_{;\bz\bz}R_{;\bz}  +  
\frac{\cplxi}{2}g^2\bar A_{;\bz}R_{;\bz\bz}  + 
 \frac{1}{4}K_{;z}g\Delta R_{;\bz}\;+
\qquad\qquad\qquad\qquad&\\\textstyle
(\frac{\cplxi}{2}g^2 A_{;zz} + \frac{1}{2}gRK_{;\bz} ) R_{;z}
+\frac{\cplxi}{2}g^2A_{;z}R_{;zz}
-K_{;\bz}(\frac{1}{4}g \Delta R_{;z} + \frac{1}{2}gRR_{;z} )
&\textstyle
= \frac{\cplxi}{2} g^2\cdot\Delta D_0.
\end{split}
\end{equation*}
Rearranging the terms and dividing by $\frac{1}{4}g$, we obtain
\begin{equation*}
K_{;z}\Delta R_{;\bz} -K_{;\bz}\Delta R_{;z}\;=\;
\textstyle
2\cplxi g \Delta D_0 -
2\cplxi g(A_{;zz}R_{;z} +A_{;z}R_{;zz}  +
\bar A_{;\bz\bz}R_{;\bz}+\bar A_{;\bz\bz}R_{;\bz}),
\end{equation*}
which finally yields
\[\textstyle
\{K, \Delta R\} = 
 \Delta D_0  - 2\Re\big((A_{;z}R_{;z})_{;z}\big).
\]
We denote this equation by $\cale^*_1$ and set
\begin{equation}\label{phi1*}\textstyle
 \varphi^*_1:=\Delta R \qquad \quad 
 D^*_1:=  \Delta D_0  - 2\Re\big((A_{;z}R_{;z})_{;z}\big). 
\end{equation}
The equation $\cale^*_1$ has the already familiar form 
\begin{equation}\label{Ersh*}\textstyle
\{K, \varphi^*_1\} = D^*_1, 
\end{equation}
and so we can apply Lemma \ref{lin-alg}. This gives us two more equations
\begin{equation}\label{Eq2*Eq3*}
\cale^*_2:=\{K,\varphi^*_2\}-D^*_2
\qquad \qquad 
\cale^*_3:=\{K,\varphi^*_3\}-D^*_3
\end{equation}
in which
\begin{align}\label{phiD2*}
\varphi^*_2&:=\{\varphi_0,\varphi^*_1\}
& 
D^*_2 &:=\{D_0,\varphi^*_1\} -\{D^*_1,\varphi_0\}
\\
\label{phiD3*}
\varphi^*_3&\textstyle
:=\half|\nabla\varphi^*_1|^2
&
D^*_3 &:=\langle\nabla D^*_1,\nabla\varphi^*_1\rangle-4\,\Re\big(A_{;z}\cdot((\varphi^*_1)_{;z})^2\big).
\end{align}
We also define the differential expressions
\begin{equation}\label{K*-i}\textstyle
\calk^*_i:=\frac{1}{\{\varphi_0,\varphi^*_1\}}\cdot \det
\pmatr{D_0\  & -(\varphi_0)_{;x^i} \\ D^*_1\ & -(\varphi^*_1)_{;x^i}}.
\end{equation}

\begin{prop}\label{prop-thm2} 
 Let $g=\lambda(x,y)(dx^2+dy^2)$ be a metric and $A$ a holomorphic
 $3$-codifferential.  Assume that $g$ satisfies the differential condition
 $\varphi_2\equiv0$ and that $\varphi^*_2=\{R,\Delta_gR\}_g$ is non-vanishing.

 Then $g$ admits a cubic integral of the form $F= \Re( A\cdot p^3 + b \cdot p^2\barp)$
 with the given tensor $A$ if and only if $\lambda$ and $A$ satisfy the covariant PDEs
 $\calg^*_3,\calg^*_2$, and $\cald$. Moreover, the component
 $b=b_1+\cplxi b_2$ is given by the formulas $b_1=\lambda^{-2}\calk^*_2,
 b_2=-\lambda^{-2}\calk^*_1$. 
\end{prop}

\proof  %
As above, the existence of such a cubic integral $F$ is equivalent to
solvability of the equation $\cale_z$:
$K_{;\bz\bz}=-\frac{\cplxi}{2}g^2A_{;z}$. Since $\cald,\calg^*_2,\calg^*_3$
are differential consequences of $\cale_z$, they are necessary conditions.

Vice versa, the system $(\cale_0,\cale^*_1)$ is solvable if and only if the
integrability condition $\calg^*_3\equiv0$ is fulfilled, and then the solution is
given by the formulas \eqref{K*-i}. Further, by Lemma \ref{lin-alg}, {\bf(c)},
the equation $\cale_z$ is a linear algebraic consequence of the conditions of
the 1st order derivatives $(\cale_0)_{;i}, (\cale^*_1)_{;i}$ and the equations
$\cale_1, \cale^*_2,\cale^*_3$. As we have shown, in the presence of
$\cale_0,\cale^*_1$ those three equations are equivalent to
$\cald,\calg^*_3,\calg^*_2$, respectively. Thus under the conditions
$\cald=0,\calg^*_2=0,\calg^*_3=0$ the solution $K$ of the system
$(\cale_0,\cale^*_1)$ solves also the equation
$K_{;\bz\bz}=-\frac{\cplxi}{2}g^2A_{;z}$. The proposition follows.
\qed

\subsection{Metrics admitting a Killing vector.}
\label{Killing}

In previous paragraphs  we considered the cases when one of the differential
expressions $\varphi_2=\{R,\half|\nabla R|^2_g\}_g$ or $\varphi_2^*=\{R,\Delta_gR\}_g$ is non-zero.  Here
we treat the problem of detecting of cubic integrals in the case when both
$\varphi_2$ and $\varphi_2^*$ vanish. Recall that by  Bonnet-Darboux-Eisenhart theorem (see
§\,\ref{previ}) this means that the metric $g$ admits a Killing vector
field $L^i\frac{\partial}{\partial x^i}$, and then  $L=L^ip_i$ is a non-trivial linear
integral. 

We maintain the notation introduced above. Define the expressions
 $\cald^*_x,\cald^*_y,\cald^*_z,\cald^*_\bz$ 
from the relations
\[\textstyle
\det
{\scriptstyle
\left(\begin{matrix}
 d\varphi_0 & D_0 \\ d\varphi^*_1 & D^*_1
\end{matrix}\right)}
= \cald^*_zdz +\cald^*_\bz d\bz
= \cald^*_xdx +\cald^*_y dy.
\]
In particular,
\begin{equation}\label{D*-eq} 
\cald^*_x := (\varphi_0)_{;x}D^*_1-(\varphi^*_1)_{;x}D_0 =
R_{;x}\cdot\big( \Delta D_0  - 2\Re((A_{;z}R_{;z})_{;z})\big)
-\Delta R_{;x}\cdot D_0
\end{equation}
and similarly for  $\cald^*_y,\cald^*_z,\cald^*_\bz$.

\begin{prop}\label{A-kill}
 Let $g$ be a Riemannian metric on a surface $S$ with the curvature $R$ and $A$
 a holomorphic $3$-codifferential. Assume that $R_{;x}$ is non-vanishing and
 that the expressions $\varphi_2=\{R,\half|\nabla R|^2_g\}_g$ and $\varphi^*_2=\{R,\Delta_gR\}_g$
 vanish. Then $A$ is compatible with $g$ if and only if it satisfies the
 equations $\cald_x=0$ and $\cald^*_x=0$.
\end{prop}

\proof  Applying the
 prolongation-projection method, we repeat the proof of Proposition
 \ref{prop-thm2} until computation of the equation $\cale^*_1$.  At this point
 the hypothesis $\varphi^*_2\neq0$ of the part {\bf(c)} of Lemma \ref{lin-alg} is not
 satisfied. The situation here is similar to that in Section \ref{degen}
  where the hypothesis $\varphi_2\neq0$ was not fulfilled.

 At this step prolongation-projection method produces the equation
 $\cald^*_x=0$, which is an analogue of the equation $\cald_x=0$ in the case
 $\varphi^*_2\equiv0$. Indeed, the condition $\varphi^*_2\equiv0$ means the degeneration of the
 matrix of coefficients of the equations $\cale_0,\cale^*_1$ whereas the pair
 of equations $(\varphi^*_2,\cald^*_x)$ means the degeneration of the
 \emph{extended} matrix of coefficients of $\cale_0,\cale^*_1$. In particular,
 under the hypotheses $R_{;x}\neq0$ and $\varphi^*_2=0$ the equations $\cale^*_1=0$ and
 $\cald^*_x=0$ are equivalent.

Since $\cald_x,\cald^*_x$ are differential consequences of the equations
$\cale_z,\cale_\bz$, their vanishing is a necessary condition for solvability
of $\cale_z$.

Vice versa, assume that $\cald_x\equiv\cald^*_x\equiv0$. Let $\cals^\#$ be the set of the
following PDEs: $\cale_z,\cale_\bz$, their derivatives up to order $2$, the
equation $\cale_0$ and its derivative $(\cale_0)_{;x}$. Then $\cals^\#$ is
involutive and hence solvable, and every solution $K$ produces a cubic
integral.

\smallskip%
The whole set of solutions of $\cals^\#$ can be constructed in the same way as
it was done in the Proposition \ref{R=const}: We fix the initial value
$K_{;x}(P)$ at some point $P\in S$, solve the ODE
$\frac{\partial}{\partial x}K_{;x}=\calk^\#_{xx}$ with the initial value $K_{;x}(P)$ along
$x$-axis, then the ODE
$\frac{\partial}{\partial y}K_{;x}=\calk^\#_{xy}$ with the initial value $K_{;x}(x,0)$ at
$y=0$, and finally express  $K_{;y}$ using  $K_{;y}=\calk^\#_y$.
\qed

\smallskip

\section{Invariant expressions.}\label{invar}

The complex calculus introduced and applied in previous sections relies on the
choice of isothermic coordinates.  In this paragraph we get rid of it: Lemma
 \ref{A-form} gives an answer on the following two questions (for a fixed
 given metric):

\begin{itemize} 
\item[(1)] Given an integral $F$ cubic in momenta, how to calculate the tensor
 $\hat A$ (the real part of the Birkhoff-Kolokoltsov form)?
\item[(2)] Given a symmetric $(3,0)-$tensor $A$, how to understand whether it
 can be the real part of a holomorphic 2-codifferential.
\end{itemize}  
 
If the coordinates we are working in are isothermal, it is easy to answer both
questions using the definition. But if the coordinates are generic, the
questions are not that trivial.
  
We would like to remark that the first question will be especially interesting
in view of program for searching new metrics admitting cubic we integrals
suggest in the conclusion Section \ref{conclu}.

We use the following notation. $x=x^1$ and $y=x^2$ are local coordinates,
$g=g_{ij}dx^idx^j$ is the metric tensor, $\lambda:=\sqrt{\det(g_{ij})}$ is the
volume density, so that $\omega_g=\lambda dx\land dy$ is the volume form, and $J_i^j$ is the
operator of the complex structure, ie., the operator of rotation by $90^\circ$ in 
the tangent bundle (it is easy to construct it for an arbitrary metric). 

\begin{lem}\label{A-form}
 \sli Let $F=(F^{ijk})$ be a symmetric $3$-vector (ie., $(3,0)$-tensor). Then
 it can be uniquely decomposed into the sum $F=\hat A+ \hat B$ where %
 $\hat A=(\hat A^{ijk})$ and $\hat B=(\hat B^{ijk})$ are also symmetric
 $(3,0)$-tensors, and $\hat A=\Re(A)$, $\hat B=\Re(B)$ for some sections $A$ of the
 bundle $T^{[3,0]}S$ and $B$ of the bundle $T^{[2,1]}S$.  Moreover, for
 arbitrary $1$-forms $\alpha_1,\alpha_2,\alpha_3$ one has
{\small\begin{align}\label{AforF} \hat A(\alpha_1,\alpha_2,\alpha_3) &=
   \frac14(F(\alpha_1,\alpha_2,\alpha_3) - F(J\alpha_1,J\alpha_2,\alpha_3)- F(J\alpha_1,\alpha_2,J\alpha_3)
   - F(\alpha_1,J\alpha_2,J\alpha_3))\\
   \hat B(\alpha_1,\alpha_2,\alpha_3) &= \frac14(3F(\alpha_1,\alpha_2,\alpha_3) + F(J\alpha_1,J\alpha_2,\alpha_3)+
   F(J\alpha_1,\alpha_2,J\alpha_3) + F(\alpha_1,J\alpha_2,J\alpha_3))
\end{align}}
or in the index form 
{\begin{align}\label{AforFindx}
\hat A^{ijk}&= \frac14(F^{ijk} - F^{i'j'k}J_{i'}^iJ_{j'}^j 
- F^{i'jk'}J_{i'}^iJ_{k'}^k - F^{ij'k'}J_{j'}^jJ_{k'}^k )
\\
\hat B^{ijk}&= \frac14(3F^{ijk} +F^{i'j'k}J_{i'}^iJ_{j'}^j 
+ F^{i'jk'}J_{i'}^iJ_{k'}^k + F^{ij'k'}J_{j'}^jJ_{k'}^k )
\end{align}}

\slii Let $A$ be a section of the bundle $T^{[3,0]}S$ and
 $\hat{A}:=\Re(A)$ its real part. Then the imaginary part $\Im(A)$ is given by
\begin{equation}\label{Im-A}
(\Im(A))^{ijk} = \frac13\bigl(\hat
A^{ljk}J_l^i + \hat A^{ilk}J_l^j + \hat A^{ijl}J_l^k).
\end{equation}

Further, $A$ is holomorphic if and only if the tensor $\hat A=(\hat A^{ijk})$
satisfies the equation $\nabla \hat A(J\cdot,J\cdot,J\cdot,J\cdot)= -\nabla \hat A(\cdot,\cdot,\cdot,\cdot)$ or in the
index form
\begin{equation}\label{A-hol-form}
A^{(ijk\,;l)}=-A^{i'j'k'\,;l'}J_{i'}^{(i}J_{j'}^jJ_{k'}^kJ_{l'}^{l)},
\end{equation}
where $(ijkl)$ means the symmetrisation in the induces. 

\sliii The principle equation $K_{;\bz\bz}=-\frac{\cplxi g^2}{2}A_{;z}$ can be
written as 
\begin{equation}\label{EqzTensor}\textstyle
\frac12 K_{;kl}( \delta^k_i\delta^l_j+J^k{}_iJ^l{}_j)+
g_{k(i}\omega_{j)l}\hat A^{klm}{}_{;m}=0
\end{equation}
where $_{(ij)}$ means the
symmetrisation in indices.
\end{lem}

\proof \sli By property (\ref{pr2}) from §\,\ref{calcul}, the bundle of
\emph{complex-valued} symmetric $3$-vectors is naturally isomorphic to the sum
$T^{[3,0]}S\oplus T^{[2,1]}S\oplus T^{[1,2]}S\oplus T^{[0,3]}S$. The complex conjugation
interchanges two inner and two outer summands. Therefore every \emph{real}
symmetric $3$-vectors has the form $F=\Re(A+B)$ for uniquely defined
sections $A$ of $T^{[3,0]}S$ and ${B}$ of $T^{[2,1]}S$. 

The formulas for $\hat A=\Re(A)$ and $\hat B=\Re(B)$ are subject of linear algebra,
therefore it is sufficient to check them for the Euclidean case
$(\rr^2,g_{Euc})$ and constant tensors ${A}=\frac{\partial}{\partial z}\cdot\frac{\partial}{\partial
 z}\cdot\frac{\partial}{\partial z}$ and $B=\frac{\partial}{\partial z}\cdot\frac{\partial}{\partial z}\cdot\frac{\partial}{\partial\bz}$. In
this case $\hat A=\Re(A)=\big( (\frac{\partial}{\partial x})^3-3(\frac{\partial}{\partial x})^2\frac{\partial}{\partial
 y}\big)$, 
$\hat B=\Re(B)=\frac{\partial}{\partial x}\big( (\frac{\partial}{\partial x})^2 + (\frac{\partial}{\partial y})^2\big)$, 
and $J(\frac{\partial}{\partial x})=\frac{\partial}{\partial y},\;J(\frac{\partial}{\partial
 y})=-\frac{\partial}{\partial x}$, and an explicit verification follows.

\medskip%
\slii Consider the covariant derivative $\nabla \hat A=\Re(\nabla A)$.  
By (\ref{pr5}) from §\,\ref{calcul}, in any complex coordinate $z$ 
we obtain $\nabla A=A_{;z}dz+A_{;\bz}d\bz$. %
Again in the coordinate $z$, multiplication with
$g\inv=\lambda\inv\frac{\partial}{\partial z}\frac{\partial}{\partial\bz}$ gives 
$g\inv\nabla A=A_{;z}\lambda\inv\frac{\partial}{\partial\bz}+ A_{;\bz}\lambda\inv\frac{\partial}{\partial z}$. Thus
$g\inv\nabla A$ is the sum of two components of the type $T^{[4,0]}S$ and
$T^{[3,1]}S$, in which the $T^{[4,0]}S$-component is ``responsible'' for
Cauchy-Riemann term $A_{;\bz}$. These two components are easily
distinguished by the operator $J$: $C(J\cdot,J\cdot,J\cdot,J\cdot)= +C(\cdot,\cdot,\cdot,\cdot)$ for
every section $C$ of the bundle $T^{[4,0]}S$, and %
$C(J\cdot,J\cdot,J\cdot,J\cdot)= -C(\cdot,\cdot,\cdot,\cdot)$ for every section $C$ of 
$T^{[3,1]}S$. Thus holomorphicity of $A$ is equivalent to the equation 
$g\inv\nabla A(J\cdot,J\cdot,J\cdot,J\cdot)= -g\inv\nabla A(\cdot,\cdot,\cdot,\cdot)$. 

Now we make use two easy observations. The first is that if $C$ is a sum of
two sections of the bundles $T^{[4,0]}S$ and $T^{[3,1]}S$ respectively, then
the $T^{[4,0]}S$-component of $C$ vanishes if and only if %
$\hat C:=\Re C$ satisfies the equation $\hat C(J\cdot,J\cdot,J\cdot,J\cdot)= -\hat
C(\cdot,\cdot,\cdot,\cdot)$. This fact can
be deduced from the action of the complex conjugation on bundles
$T^{[p,q]}S$, see (\ref{pr3}) from §\,\ref{calcul}. Another observation is
that in the index form the symmetric tensor $\Re(g\inv\nabla A)$ is given by
$A^{(ijk\,;l)}=A^{(ijk}{}_{;m}g^{l)m}$. This implies the formula
\eqref{A-hol-form}.

\smallskip%
Now consider the tensor $A_{;z}dz$. By the construction, this is the
section of the bundle $T^{[3,0]}S\otimes\Omega^{[1,0]}S$. The latter bundle is isomorphic
to $T^{[2,0]}S$, and the isomorphism is induced by the isomorphism
$T^{[1,0]}S\otimes\Omega^{[1,0]}S\cong\cc$. The latter isomorphism is simply the operation of
contraction of indices. This operation commutes with complex conjugation, and
therefore with the operator $\Re$ of taking real part.  It follows that the
image of $\Re(A_{;z}dz)$ in $T^{[2,0]}S$ is given by the symmetric
$(2,0)$-tensor $\Div{\hat A}$ with components $(\hat A^{ijk}{}_{;k})$. 
It should be noticed
that besides the usual symmetry $\hat A^{ijk}{}_{;k}=\hat A^{jik}{}_{;k}$ 
this tensor has another symmetric property, 
namely $\Div{\hat A}(J\cdot,J\cdot)= -\Div{\hat A}(\cdot,\cdot)$. Thus $\Div{\hat A}$
has two independent components, as it should be since $A_{;z}$ is complex-valued.

The last formula expressing $\Im(A)$ can be deduced in a similar way as above
and is left to the reader as an easy exercise (and is not important for us).

\medskip%
\sliii The real part of the principle equation \eqref{Eqz} is 
\begin{equation}\label{ReEqz}\textstyle
K_{;\bz\bz}+K_{;zz}= \frac{\cplxi g^2}{2}(\bar A_{;\bz}-A_{;z}).
\end{equation}
Since the principle equation  is the $\Omega^{[0,2]}$-component of
\eqref{ReEqz}, both equations are equivalent. Since $T(Jv,Jw)= T(v,w)$ for
every section $T$ of the bundle $\Omega^{[1,1]}S$ and $T(Jv,Jw)= -T(v,w)$ for
sections of $\Omega^{[2,0]}S\oplus\Omega^{[0,2]}S$, the tensor $K_{;\bz\bz}+K_{;zz}$ has
components 
\[\textstyle
\frac12 K_{;kl}( \delta^k_i\delta^l_j+J^k{}_iJ^l{}_j).
\]

Rewriting the right hand side of \eqref{ReEqz} we start with the equality
$\frac{\cplxi}{2}(\bar A_{;\bz}-A_{;z}) = \Im(A_{;z})$.  Since $A$ is
holomorphic and $A_{;\bz}=0$, the tensor $A_{;z}$ equals
$A^{ijk}{}_{;k}$. Using the fact that the isomorphism \eqref{trace} is simply
the contraction of indices, we obtain the formula
$g_{ik}g_{jl}\Im(A^{klm}{}_{;m})$.  Further, since $A_{;z}$ has type
$T^{[2,0]}S$, $\cplxi{\cdot}A_{;z}=J{\cdot}A_{;z}$ and hence
$\Im(A_{;z})=-\Re(J{\cdot}A_{;z})$. Using the equality $g_{ij}J^j{}_k=\omega_{ik}$ we
obtain the desired form \eqref{EqzTensor} of the equation.
\qed

\medskip%

\section{Pseudo-Riemannian case.}\label{sign=1,1}

In this section we show that the results and formulas obtained in the case of
Riemannian metrics remain valid in the pseudo-Riemannian case after an 
appropriate modifications which we describe.

\subsection{Null-coordinates.}  These are the counterpart of isothermic coordinates
in the pseudo-Riemannian case.

\begin{defi}\label{null-coord} \rm Let $g$ be a pseudo-Riemannian metric on a
 surface. A vector (field) $v$ satisfying $g(v,v)=0$ is called a
 \slsf{null-vector (field)}. \slsf{Null-coordinates}%
 \footnote{From the mathematical point of view, \slsf{isotropic coordinate
   system} would be probably a better notion. However, this terminology is
  already established in the physical literature, see e.g.\ \cite{Cu-La}, and
  moreover, the notion \slsf{isotropic coordinates} itself has another
  meaning, see e.g.\ \cite{Cro}.}
 are such coordinates in which $g$ has the has anti-diagonal form
 $g=\lambda(x,y)dxdy$.
\end{defi}

\begin{lem}[Folklore]\label{coord-1,1} Let $g$ be a pseudo-Riemannian metric on a
 surface $S$. Then at each point on $S$ there exist local coordinates
 $x=x^1,y=x^2$ in which $g$ has the anti-diagonal form $g=\lambda(x,y)dxdy$
\end{lem} 

\proof At each point $p$ on $S$ the metric $g$ has exactly two isotropic
directions, ie., pointwisely linearly independent vectors fields $v,w$ with
$g(v,v)=g(w,w)=0$. Integrating them we obtain a $2$-web of curves on
$S$. The local parameters functions for this $2$-web are the desired
coordinates $x=x^1,y=x^2$.
\qed

\medskip%
Let us locally fix the order and the orientations of null-directions.
Then an null-coordinate system is defined uniquely up to oriented
reparametrisations $x'{}^i=f^i(x^i)$.

\subsection{Calculus in null-coordinates.} \label{null-calcul} As in the
Riemannian case, formulas and calculus become much simpler if we use
null-coordinates.

In analogy with the complex calculus, define \emph{real} bundles $T^{[1,0]}S$,
$T^{[0,1]}S$, $\Omega^{[1,0]}S$, $\Omega^{[0,1]}S$, generated respectively by vectors
$\frac{\partial}{\partial x^1},\frac{\partial}{\partial x^2}$ and the forms $dx^1$, $dx^2$. The definition
of these bundles depends only on the local combinatorial data, whereas the
decompositions $TS=T^{[1,0]}S\oplus T^{[0,1]}S$ and $T^*S=\Omega^{[1,0]}S\oplus\Omega^{[0,1]}S$ are
defined solely by the metric $g$. Further, for $p,q\geq0$ set
\[
T^{[p,q]}S:=(T^{[1,0]}S)^{\otimes p}\otimes(T^{[0,1]}S)^{\otimes q}
\qquad\text{and}
\qquad
\Omega^{[p,q]}S:=(\Omega^{[1,0]}S)^{\otimes p}\otimes(\Omega^{[0,1]}S)^{\otimes q}.
\]
Then 
\begin{equation}\label{Tpq}
T^{[p,q]}S \otimes T^{[p',q']}S =T^{[p+p',q+q']}S
\end{equation}
for non-negative $p,p',q,q'$. We extend the definition of the bundles
$T^{[p,q]}S$, $\Omega^{[p,q]}S$ for all integers $p,q$ by means of the formula
\eqref{Tpq} and the relation
\begin{equation}\label{Ohmpq}
\Omega^{[p,q]}S=T^{[-p,-q]}S.
\end{equation}
Every section of the bundle $T^{[p,q]}S$ or $\Omega^{[p,q]}S$ is given by its
single component/coefficient which is a real function.

As in the case of complex calculus, we imbed the bundles $T^{[p,q]}S$,
$\Omega^{[p,q]}S$ with $p,q\geq0$ in the bundles of \emph{symmetric} $(k,0)$- and
$(0,k)$ tensors. The imbedding is done using symmetrisation operator. This
gives us natural decompositions $\sym^kTS=\oplus_{p+q=k}^kT^{[p,q]}S$ and
$\sym^kT^*S=\oplus_{p+q=k}^k\Omega^{[p,q]}S$.

The metric $g$ is a section of the bundle $\Omega^{[1,1]}S$. We use multiplication
with powers of $g$ to define natural isomorphisms
$T^{[p,q]}S\xrar{\;\cong\;}T^{[p+k,q+k]}S$ and
$\Omega^{[p,q]}S\xrar{\;\cong\;}\Omega^{[p+k,q+k]}S$.

Further, define the operators $\nabla^{[1,0]}:\Omega^{[p,q]}S\to:\Omega^{[p+1,q]}S$ and
$\nabla^{[0,1]}:\Omega^{[p,q]}S\to:\Omega^{[p,q+1]}$ as the corresponding homogeneous
components of the operator $\nabla$. Thus for any section $F=f(x,y)dx^pdy^q$ of the
bundle $\Omega^{[p,q]}S$ with $p,q\geq0$ we obtain 
\[
\nabla^{[1,0]}F= f_{;x}dx^{p+1}dy^q,
\qquad\qquad
\nabla^{[0,1]}F= f_{;y}dx^{p}dy^{q+1},
\]
and $\nabla^{[1,0]}F+\nabla^{[0,1]}F=\sym(\nabla F)$ where $\sym$ denotes the symmetrisation
operator. In the case of a section $A=a(x,y)(\partial_x)^p(\partial_y)^q$ of the
bundle $T^{[p,q]}S$ with $p,q\geq0$ we obtain 
\[
\nabla^{[1,0]}A= a_{;x}(\partial_x)^{p+1}(\partial_y)^q 
\qquad\qquad
\nabla^{[0,1]}A= a_{;y}(\partial_x)^p(\partial_y)^{q+1},
\]
and $\nabla^{[1,0]}A+\nabla^{[0,1]}A=\tr(\nabla A)=\Div(A)$ where $\tr$ denotes the
contraction of indices. The operators $\nabla^{[1,0]}$ and $\nabla^{[0,1]}$ do not
commute, and the commutator formula is similar to the complex case
\eqref{commutator}:
\begin{equation}\label{commut-psR}\textstyle
\begin{split}
 (\nabla^{(1,0)}\nabla^{(0,1)} - \nabla^{(0,1)}\nabla^{(1,0)})fdx^pdy^q &=
 \frac{q-p}2\cdot R\cdot g\cdot fdx^pdy^q\\
 (\nabla^{(1,0)}\nabla^{(0,1)} - \nabla^{(0,1)}\nabla^{(1,0)})f(\partial_x)^p(\partial_y)^q
&= \frac{p-q}2\cdot R\cdot g\cdot  f(\partial_x)^p(\partial_y)^q
\end{split}
\end{equation}

\smallskip%
Let $F(x,y;p_x,p_y)$ be a function on $T^*S$ cubic in momenta. Then it has a
form $F=A_1p_x^3-A_2p_y^3+(B_1p_x+B_2p_y)p_xp_y$ for uniquely defined sections
$A_1$ of $T^{[3,0]}S$, $A_2$ of $T^{[0,3]}S$, $B_1$ of $T^{[2,1]}S$, and $B_2$
of $T^{[2,1]}S$. The Hamiltonian function corresponding to the metric has the
form $H=\frac{p_xp_y}{2\lambda(x,y)}$, and the equation of integrability of $F$
reads 
{\small\begin{equation}\label{braFH-psR}
\begin{split}
\textstyle
\{F,H\} &\textstyle
= \frac{1}{2\lambda^2}\Big(\,
p^4_x\, \lambda \cdot (A_1)_{;y} - p^4_y \, \lambda \cdot (A_2)_{;x}
+ 
\\[3pt]
&
p^3_xp_y \cdot ( B_1 \cdot \lambda_{;y} +3A_1 \cdot \lambda_{;x}+ \lambda \cdot  (B_1)_{;y}  +\lambda  \cdot (A_1)_{;x})+\\[3pt]
&
p_xp^3_y \cdot ( B_2 \cdot \lambda_{;x} -3A_2 \cdot \lambda_{;y}+ \lambda  \cdot (B_2)_{;x}  -\lambda  \cdot (A_2)_{;y}) 
+\\[3pt] &
p_x^2p_y^2  \cdot (2B_1 \cdot \lambda_{;x} + \lambda \cdot (B_1)_{;x} +  2B_2 \cdot \lambda_{;y} + \lambda  \cdot (B_2)_{;y} )
\Big)
\end{split}
\end{equation}}

\begin{lem}\label{kolokol}  Let $g=\lambda dxdy$ be a pseudo-Riemannian
 metric on a surface $S$ in null-coordinates.

\smallskip\noindent %
 {\bf(a)} A function $a(x,y)$ is independent of the variable $y$ if and
 only if for some/every $k\geq0$ the tensor $A:=a\,\partial_x^k$ satisfies the equation
 $A_{;y}=0$.  

 \smallskip\noindent %
 {\bf(b)} Let $F$ be an integral of $g$ cubic in momenta. Then its
 components $A_1$ and $A_2$ depend only on one corresponding variable.

 Moreover, if $A_1$, $A_2$ are non-zero at a given point on $S$, then there
 exists a null-coordinate system $(x,y)$, unique up to exchange of the
 coordinates, in which components of $A_i$ are identically $1$, ie.,
 $A_1\equiv\partial_x^3$ and $A_2\equiv\partial_y^3$.

\end{lem} 

Here $A_{;y}$ means the covariant derivative $\nabla^{[0,1]}A$.

\begin{rem} Similar statement  holds for  polynomial integrals of  other degree, see \cite[§3.3]{Bo-Ma-Pu1} for the case of quadratic integrals. 
\end{rem}  
\proof \ If $(x',y')$ is any null-coordinate system with
$A_1=a_1(x')\partial_{x'}^3$ and $A_2=a_2(y')\partial_{y'}^3$, then the desired coordinates
are defined by $x=\int\frac{dx'}{\sqrt[3]{a_1(x')}}$ and
 $y=\int\frac{dy'}{\sqrt[3]{a_2(y')}}$.
\qed

\begin{defi}\label{def-kolo} \rm For any function  $F$ on the cotangent bundle
 $T^*S$ cubic in momenta we call the components $A_1$ and $A_2$ the {\it
  Birkhoff-Kolokoltsov} $3$-codifferentials associated with $F$ and
 pseudo-Riemannian metric $g$.
 $A_1$ and/or $A_2$ is said to be {\it
  quasi-holomorphic}%
 \footnote{In view of {\it pseudo}-Riemannian metrics, {\it pseudo}-holomorphic
  would be probably a better notion. However, this terminology is already
  reserved, see \cite{Gro}.}  %
if $A_i=a_i(x^i)\partial_{x^i}^3$ with $a_i(x^i)$ depending only on one variable.

The null-coordinate system $(x,y)$ in which
 $A_1\equiv\partial_x^3$ and $A_2\equiv\partial_y^3$ is called \slsf{adapted} to $A_1,A_2$ or to $F$.
\end{defi}

\begin{prop}\label{Eq-psR} Let $g=\lambda dxdy$ be a pseudo-Riemannian
 metric on a surface $S$ in null-coordinates and $F$ a function cubic in
 momenta. Then $F$ is an integral for $g$ if and only if 
\begin{itemize}
\item[(A)] its components $A_1,A_2$ are quasi-holomorphic;
\item[(B)] there exists a function $K(x,y)$ such that  the components
 $B_1,B_2$ are related  as $B_1=-\lambda^{-2}K_{;y}$,  $B_2=\lambda^{-2}K_{;x}$; 
\item[(K)] the function $K$ satisfies the equations 
\[
K_{;xx}=g^2\,(A_2)_{;y} \qquad\qquad K_{;yy}=g^2\,(A_1)_{;x}.
\]
\end{itemize} 
\end{prop}

\proof Five homogeneous components of bi-degree $T^{[i,j]}S$ ($i+j=4$) of the
Poisson bracket \eqref{braFH-psR} are five equations on $A_i,B_i$ and $\lambda$. The
function $F$ is a cubic integral for $g=\lambda dxdy$ iff all five equations are
satisfied. The first two are equations (A). The last one can be written as 
$d(\lambda^2B_1dx + \lambda^2B_2dy)=0$ and hence is equivalent to the  property
(B). Finally, substituting expressions (B) into the remaining two equations 
(second and third lines in \eqref{braFH-psR}) we obtain equations (K).
\qed

\subsection{Calculation of equations  in pseudo-Riemannian case.}
Here we use the trick and translate all the results of 
§\ref{prol1} and Section \ref{degen} almost literally in the new situation. Doing so,
we must follow the next rules.
\begin{itemize}
\item the complex coordinate $z$ is replaced by the null-coordinate $x$,
 the coordinate $\bz$ by $y$;
\item the holomorphic 3-codifferential $A$ is replaced by the
 quasi-holomorphic 3-codifferential $A_1$, $\bar A$ by $A_2$,
\item the imaginary unit $\cplxi$ is dropped (replaced by $1$), the operator
 $J^i{}_j$ of rotation by $90^\circ$ is replaced by the identity $\id$ operator;

\item the complex conjugation $a\mapsto\bar a$ is replaced by the
 \slsf{null-involution} $I:(x,y)\mapsto(y,x)$; the operator
 $\Im(A)=\frac{1}{2\cplxi}(a-\bar a)$ transforms into the operator
 $P_-:a\mapsto\half(a- I(a))$, this is the operator of the projection on the
 eigenspace $\mathbb{E}_I(-1)$ of the operator $I$ corresponding to the
 eigenvalue $-1$; similarly, the operator $\Re(A)=\frac{1}{2}(a+\bar a)$
 transforms into the operator $P_+:a\mapsto\half(a+ I(a))$ corresponding to the
 eigenvalue $+1$;
\item the Poisson bracket $\{f,g\}=\frac{1}{2\cplxi
  g}(f_{;z}g_{;\bz}-f_{;\bz}g_{;z})=g\inv\Im(f_{;z}g_{;\bz})$ is now given by
 the new formula $\{f,g\}:=g\inv
 P_-(f_{;x}g_{;y})=\frac{1}{2g}(f_{;x}g_{;y}-f_{;y}g_{;x})$;
\item if it is possible, every equation must be written in complex form and
 then its conjugate must included to the system; exceptions are the equations
 which are real ``in their nature'', such as
 $\cale_0,\ldots,\cale_3;\cale^*_1,\ldots,\cale^*_3,\allowbreak\calg_0,\ldots,\calg_3$; then
 after applying the above rules every pair of complex conjugate equations
 transforms into a pair of equations $\cale, I(\cale)$ where $I$ is the
 null-involution above
\end{itemize} 

For example, the equation $K_{;\bz\bz}=A_{;z}$ is complemented by its
conjugate $K_{;zz}=A_{;\bz}$, so that counterpart of $K_{;\bz\bz}=A_{;z}$ in
the pseudo-Riemannian case is the pair of the equations $K_{;xx}=A_{2;y}$ and
$K_{;yy}=A_{1;x}$. In particular, we have two real equations in both
Riemannian and pseudo-Riemannian cases.

\begin{prop}\label{algo-psR} {\bf (a)} After changes made by the rules described
 above, Theorems \ref{th:main1} and \ref{th:main2} remain valid in the
 pseudo-Riemannian case.

 \smallskip\noindent%
 {\bf(b)} Let $g$ be a pseudo-Riemannian metric with vanishing $\varphi_2,\varphi^*_2$ and
 non-vanishing $\varphi_1$. Then g admits a nontrivial linear integral.
\end{prop}

\proof \;{\bf(a)} \;The proof of Theorems \ref{th:main1} and \ref{th:main2} is
done by means of formal calculus in which the condition of positivity of the
metric was never used.

\smallskip\noindent%
{\bf(b)} \; Considering in \cite[pp.\;323-325]{Ei} Riemannian metrics $g$ on
surfaces with non-constant curvature $R$, Eisenhart introduces the following
coordinate system: One of the coordinates is the curvature $R$ itself and the
other is orthogonal to the first one. In other words,
$g=g_{11}dx^2+g_{22}dy^2$ and $R(x,y)=x$. The necessary and sufficient
condition for the existence of such coordinates is that $\nabla R$ is not orthogonal
to itself. In the Riemannian case this simply means the non-vanishing of $\nabla R$
and thus the non-constancy of $R$, in the pseudo-Riemannian case this is the
condition  $\varphi_1\neq0$. The rest of Eisenhart's proof works also in pseudo-Riemannian.
\qed

\subsection{Pseudo-Riemannian metrics with  $|\nabla R|_g^2=0$: local classification and nonexistence of cubic integrals }\label{phi1=0}

The goal of this paragraph  is to prove

\begin{prop}\label{qub-phi1} If a pseudo-Riemannian metric $g$ with
 $|\nabla R|_g^2\equiv0$ admits a non-trivial cubic integral, then $R$ is constant. 
\end{prop}

We will  need the following 

\begin{lem}\label{lem-dR} Let $g$ be a pseudo-Riemannian metrics with
 $\nabla R\neq0$ and  $|\nabla R|_g^2=0$. Then in generic null-coordinates it has the form
\begin{equation}\label{g-dR}
g=\frac{f_2(x)f_1(y)_y\;dxdy}{(f_1(y)+f(x))^2}
\end{equation}
with some functions $f_1(y), f(x),f_2(x)$, and then the curvature is given by
\[
R=\frac{2f(x)_x}{f_2(x)}.
\]

Moreover, there exist
null-coordinates in which $g$ has the normal form
\begin{equation}\label{g-dR-N}
g=\frac{dxdy}{(y+f(x))^2}
\qquad\text{with}
\qquad
R=2f(x)_x.
\end{equation}
where $f(x)$ is some function.

Such a metric $g$ admits no Killing vector field, ie., no linear integral $L$.
\end{lem}

\state Remark. The generic formula is given up to exchange of null-coordinates
$(x,y)$. In \eqref{g-dR} and in the proof subscript indices
$f_1(y)_y,u_{xyy},\ldots$ denote usual (partial) and \emph{not covariant} derivatives of
functions $f_1(y),u(x,y),\ldots$

\proof Choose some null-coordinates $(x,y)$. The condition $|\nabla R|_g^2=0$ means
that the gradient $\nabla R$ is a null-vector. Hence $R$ depends only on one of
null-coordinates, say on $x$, and $R_y=0$. On the other hand, $R_x\neq0$ by the
hypotheses of the lemma. Inverting the coordinate $x$ if needed, we may assume
that $g=e^{u(x,y)}\,dx\,dy$ for some function $u(x,y)$. Then $R=e^{-u}u_{xy}$
and the equation $|\nabla R|^2_g=0$ reads $\big(e^{-u}u_{xy}\big)_y=0$, or
$u_{xyy}-u_{xy}u_y=0$. Substitution $u_y=v$ yields $v_{xy}-vv_x=0$, which is
equivalent to $(2v_y-v^2)_x=0$. Hence $2v_y=v^2+\psi(y)$ for some function
$\psi(y)$. This is a Riccati ODE, whose generic solution can be found as
follows. Let $v_1(y)$ be a fixed special solution of the equation
$2v_y=v^2+\psi(y)$, so that $\psi(y)=2(v_1(y))_y-v_1(y)^2$. Since $\psi(y)$ was
arbitrary, we may assume that $v_1(y)$ is an arbitrary function and $\psi(y)$ is
given by the above formula.  Now make the substitution
$v(x,y)=v_1(y)+\frac{2}{w(x,y)}$. Then the equation $2v_y=v^2+\psi(y)$ transforms
into $w_y+v_1w+1=0$. Since $v_1(y)$ is an arbitrary function we may assume it
has the form $v_1(y)=\frac{f_1(y)_{yy}}{f_1(y)_y}$ with an arbitrary function
$f_1(y)$.  Then the equation $w_y+v_1w+1=0$ reads
$f_1(y)_yw_y+f_1(y)_{yy}w+f_1(y)_y=0$, which means $(f_1(y)_yw(x,y)+f_1(y))_y=0$.
Consequently, $w(x,y)=-\frac{f_1(y)+f(x)}{f_1(y)_y}$ with an arbitrary function
$f(x)$, and 
\[
\textstyle
v(x,y)=  \frac{f_1(y)_{yy}}{f_1(y)_y}  -2\frac{f_1(y)_y}{f_1(y)+f(x)}.
\]
Now integration yields $u(x,y)=\log(f_2(x))+\int v(x,y)dy$ which gives the desired
form of the metric $e^{u(x,y)}dxdy$.

\smallskip
The functions $f_1(y)$, $f_2(x)$ can be easily eliminated using the substitution
$y'=f_1(y)$, $x'=\int\!\!f_2(x)dx$. 

\smallskip
Now assume that a metric $g=\lambda(x,y)dxdy$ with $R_y=0$ admits a linear
integral $L$. This is equivalent to the existence of Killing vector field $\xi$.
Since the curvature $R$ is constant on lines $x=c$,   $\xi$ must have the form
$\xi=\alpha(x,y)\partial_y$. The Lie derivation gives
\[
\scrl_\xi(\lambda\,dx\,dy)= \alpha\,\partial_y\lambda\,dx\,dy + 
\lambda\,\partial_x\alpha\,dx\,dx+\lambda\,\partial_y\alpha\,dx\,dy.
\]
In particular, $\partial_x\alpha$ must vanish. Then $\alpha$ depends only on $y$, and after the
appropriate reparametrisation $\xi=\partial_y$. But then $\lambda$ must be independent of
$x$ which would imply the vanishing of $R$.
\qed

{\it Proof of Proposition \ref{qub-phi1}.} 
We use the following notation: $x,y$ are null-coordinates in
which the metric has the form \eqref{g-dR-N}, $f=f(x),f_1(x)=f_1,f_2(x)=f_2$
are functions of the variable $x$, $f_1',f'''_2,f''$ their
derivatives. Further, we set $A_1=a_1(x)\partial_x^3$, $A_2=a_2(y)\partial_y^3$, and then
$a_1',a_2''$ denote the derivatives (in $y$ in the case $a_2$).  The usual
(non-covariant!) partial derivatives are denoted by $\partial^3_{xxy}K$ and so
on. Besides we introduce the function $Y(x,y):=y+f(x)$. Clearly, it satisfies
the relations $\partial_xY=f'$ and $\partial_yY=1$. Finally, $c_1,c_2,\ldots$ will be constants.

In the case $g=\frac{dx\,dy}{Y^2}$ with $Y=y+f(x)$ the only non-zero
Christoffel symbols are $\Gamma_{11}^1=-\frac{2f'}{Y}$ and
$\Gamma_{22}^2=-\frac{2}{Y}$. This gives us the expansions
\[
K_{;yy}= \partial^2_{yy}K+\frac{2}{Y}\partial_yK
\qquad\text{and}
\qquad
g^2A_{1\,;x}=\frac{Ya'_1-6f'a_1}{Y^5}.
\]
Thus $K_{;yy}-g^2A_{1\,;x}$ appears to be an ODE 
\[
\partial^2_{yy}K+\frac{2}{y+f}\partial_yK -\frac{(y+f)a'_1-6f'a_1}{(y+f)^5}=0
\]
on $K$ with respect to $y$ with rational coefficients in which $x$ is a
parameter. The direct integration gives
\begin{equation}\label{K-form1}
\textstyle
K(x,y)=f_1(x) +\frac{f_2(x)}{Y}  + \frac{Ya_1'-2f'a_1}{2Y^3},
\end{equation}
and hence
\begin{equation}\label{Kx-form1}
\textstyle
\partial_xK(x,y)=f'_1 +\frac{f'_2}{Y} -\frac{f'f_2}{Y^2}  
+ \frac{a_1''}{2Y^2}-
\frac{f'a'_1}{Y^3}-
\frac{f'a'_1+f''a_1\vph}{Y^3}+\frac{3f'^2a_1}{Y^4}
\end{equation}
In the same way we obtain 
\[
\textstyle
K_{;xx}= \big(\partial_x+\frac{2f'}{Y}\big)\partial_xK
\qquad\text{and}
\qquad
g^2A_{2\,;y}=\frac{Ya'_2-6a_2}{Y^5}.
\]

Now we make the following observation: Substituting of \eqref{K-form1} in
$Y^5\big(\partial^2_{xx}K+\frac{2f'}{Y}\partial_xK\big)$ and using $\partial_xY=f',\partial_yY=1$ we
obtain a polynomial (say $P(Y)$) of degree $5$ in $Y$ whose coefficients are
smooth functions in $x$, such that $a_2(y)$ satisfy the equation
$Y\partial_ya_2-6a_2=P(Y)$. For any fixed $x$ this is a polynomial ODE on
$a_2$. Integrating it, we obtain that (for any fixed $x$!) is given by
\begin{equation}\label{a2Y}
\textstyle
a_2(y)=Y^6\cdot\big( f_3(x) + \int Y^{-7}P(Y)dY\big).
\end{equation}
In view of this formula, the problem of compatibility of the
$3$-codifferential $A=(A1,A2)$ with the metric $g=\frac{dx\,dy}{(y+f(x))^2}$
is equivalent to finding the functions $f(x),f_1(x),f_2(x),f_3(x),a(x)$ whose
substitution in \eqref{a2Y} gives a function independent of $x$.

Further, if we write $P(Y)=\sum_{i=0}^5p_iY^i$, then the integration of
\eqref{a2Y} gives a polynomial expression for $a_2$ 
\[
\textstyle 
a_2=f_3Y^6 -\sum_{i=0}^5\frac{p_i}{6-i}Y^i
\]
with coefficients depending only on $x$. After expanding $Y=y+f(x)$ we still
obtain an expression for $a_2$ which is a polynomial
$Q(y)=\sum_{j=0}^6q_j(x)y^j$ in $y$ of degree $6$ with coefficients $q_j(x)$
depending only on $x$. However, since $a_2$ is independent of $y$, all these
coefficients must be constant. In this way we obtain $6$ equations
$q_i(x)=c_i$ (ODE in general) on the above functions
$f(x),f_1(x),f_2(x),f_3(x)$ and $a_1(x)$. The expressions $q_i(x)$ can be
calculated (or checked) using the formulas above, and we simply provide the
final results.

We are going to solve these equations using the prolongation(few times) and
projection(mostly) method, resolving subsequently equations
$q_6(x)=c_6,q_5(x)=c_5,\ldots$ and substituting the results into the next
equation. The obtained new ODEs will be denoted by $E_6,E_5,\ldots$ 

The first equation is simply $f_3(x)=c_6$, so $f_3$ is a constant. The next
equation $q_5(x)=c_5$ (up to constant factor) reads as  $c_5+6c_6f-f_1''=0$ or
\begin{equation}\label{f1xx}
f_1''=c_5+6c_6f.
\end{equation} 
Substituting  the latter expression  $q_4(x)=c_4$ we obtain
the next ODE  $E_4$ which reads $c_4 -15 c_6 f^2
- f_1' f'  - \half f_2'' -5 f c_5$. %
Resolving it we obtain 
\begin{equation}\label{f2xx}
f_2''=-2f_1'f'-10fc_5-30c_6f^2+2c_4.
\end{equation} 
Next one substitution gives the equation $E_3$ which can be resolved as follows:
\begin{equation}\label{a1xxx}\textstyle
a_1''' = 6 c_3   +120 c_6 f^3
+ 2 f_2 f'' +60 f^2 c_5   -24 f c_4.
\end{equation}
Doing the same with $q_2(x)=c_2$, we obtain 
\begin{equation}\label{a1xx}
a_1'' = (2f')\inv \big(
 -2 c_2   -3 a'_1 f'' -a_1 f'''
+12 f c_3   +60 c_6 f^4 + 40 f^3 c_5- 24 f^2 c_4 
\big).
\end{equation}
The result of the same procedure for $q_1(x)=c_1,q_0(x)=c_0$ is
\begin{equation}\label{a1x}
a_1' = (10f'^2)\inv \big(
30 f ^2 c_3   +60 c_6 f^5   -
14 a_1 f''  f'- 40 f^3 c_4+50 f^4 c_5+5 c_1-10 f c_2 
\big),
\end{equation}
\begin{equation}\label{a1}
a_1 = (2f'^3)\inv \big(
2 f^4 c_4-2 f^3 c_3+c_0-2 f^5 c_5 
-2 c_6 f^6-f(x) c_1+f^2 c_2
\big).
\end{equation}

Recall that we have substituted in each successive equation the results of
preceding calculations. This was clearly the projection procedure, and the
next one is the prolongation one. Let us denote the expressions in
(\ref{a1xxx}--\ref{a1}) by $\cala_3,\ldots,\cala_0$, so that the formulas read
$a_1^{(i)}=\cala_i$.  We compute solely the consistency equation
$\cala_1-\partial_x\cala_0$. The result is
\begin{equation}\label{cons-a1x}
\textstyle
\cala_1-\partial_x\cala_0= \frac{f''}{10 f'^4}\big(
16 c_6 f^6  +16 f^5 c_5
 -16 f^4 c_4+16 f^3 c_3-8 f^2 c_2
+8 f c_1  -8 c_0.
\big)
\end{equation} 
Since this consistency condition must be satisfied identically, we conclude
that there are two possibilities: either $f''$ vanishes identically, or all
constants $c_0,\ldots,c_6$ must be zero. By \eqref{g-dR-N}, $f''\equiv0$ means that the
curvature $R$ is constant. 

We exclude this possibility and consider the alternative case $c_0=\ldots=c_6=0$.
Then from \eqref{a1} we see that $a_1$ vanishes identically.  The equation
\eqref{a1xxx} reads now $f_2f''=0$, and so $f_2$ also vanishes.  Finally, from
\eqref{f2xx} we obtain $f_1'f'=0$, which means that $f_1(x)$ is constant.

Summing up, we conclude from \eqref{K-form1} that the the function $K(x,y)$ is
constant. Finally, let us observe that the equation
$(y+f(x))\partial_ya_2(y)=6a_2(y)$ admit no non-trivial solution depending only on
$y$. \qed

\section{Conclusion}\label{conclu}

We presented an algorithm that, given a pair $(g, A)$, where $g $ is a
two-dimensional metric, and $A $ is a holomorphic 3-codifferential, answers
the question whether there exists a cubic integral whose Birkhoff-Kolokoltsov
3-codifferential coincides with $A$. Moreover,  in the most interesting
 cases covered by Theorems \ref{th:main1}, \ref{th:main2},  it provides
precise formulas for the integral. The algorithm works in an arbitrary
coordinate system: we need to calculate certain expressions given by precise
algebraic formulas including the components of $g$ and $A$ and their
derivatives, and compare them with zero.  It is easy to implement the
algorithm on modern computer algebra packages, say Maple\textsuperscript{®} or
Mathematica\textsuperscript{®}.

Our results suggest the following program for search for new natural
Hamiltonian system admitting integrals polynomial in momenta integrals of
degree $\leq 3$.

Let us take a metric on a surface whose geodesic flow admits a cubic integral.
For example, we can take a metric of constant curvature, or  a metric
admitting a Killing vector field, or a metric coming from one of the known
natural Hamiltonian system via Maupertuis principle.

Let $K:=|\vec p|_g$ be the kinetic energy corresponding to such metric and
$F_3$ be the integral. Let us now look for a function $F_1:T^*S\to \mathbb{R}$
linear in momenta, and for a function $V:S\to \mathbb{R}$ such that $F_3+F_1$ in
an integral for $K + V$, i.e., $\{F_3+F_1, K+ V\}=0$, where $K$ is the kinetic
energy corresponding to $g$.  By Maupertuis principle, the later condition is
equivalent to the statement that for every constant $h$ the functions
$H_h:=\tfrac{1}{(h-V)} K_g$ and $F_h:= F_3 + H_h\cdot F_1$ commute: $\{H_h,
F_h\}=0$.

Clearly, the Birkhoff-Kolokoltsov 3-codifferential of the pair $H_h, F_h$ does
not depend on the parameter $h$ and coincides with the Birkhoff-Kolokoltsov
3-codifferential of the pair $K, F_3$, i.e., is known since we know $K$ and
$F_3$. Thus, the conditions $\calg_2=0$, $\calg_3=0$ (or $\mathcal{D}=0$,
$\calg^*_2=0$, $\calg^*_3=0$ in the case covered by Theorem \ref{th:main2}) can be viewed as quasilinear PDEs on the
coefficients of $F_1$ and on $V$, i.e., on three unknown functions. Since the
conditions are fulfilled for every $h$, the system of PDEs for $V$ and $F_1$
is overdetermined. Easy analysis shows that it is of finite type and there
exists an algorithmic (though highly computational) way to find a solution.
  
In other words, we suggest to look for new systems by adding potential to the
known integrable geodesic flows, and the results of our paper ensures that one
can do it algorithmically.
    
Note that all known natural Hamiltonian systems either come from physics (e.g.
Goryachev-Chaplygin), or were obtained by the naive version of the above
procedure (for example those obtained in \cite{Du-Ma,Se,Ve-Ts,Ye}).
  
As the most promising metrics $g$ to start the above program, we consider the
metrics admitting linear integrals. Recall that such metric are completely
described. Local description is due at least to \cite[§§590--593]{Da}, global
(= when the manifold is closed) can be found for example in \cite{Bo-Ma-Fo}.
Note that the starting point for the systems from \cite{Du-Ma,Se,Ve-Ts,Ye}
were certain metrics admitting linear integrals.  This will be the next
direction of our research; we are quite positive that it is possible to
describe all natural system admitting integrals polynomial in momenta of
degree $\leq 3$ such that the corresponding metric admits a Killing vector field,
and hope to find new examples.
  
Another perspective application is to try to prove/disprove the following
conjecture from \cite{Bo-Ko-Fo}: {\it if a real-analytic metric on the 2-torus
 admits a cubic integral, then it admits a linear integral,} or its weaker
form: {\it if a real-analytic natural Hamiltonian system on the 2-torus admits
 an integral of degree three in momenta, then it admits a linear integral. }
Indeed, as we explained in Corollary \ref{g-geq-2}(ii), the space of
holomorphic 3-codifferential on the torus is a 2-dimensional linear
space. More precisely, there exists a global (periodic) coordinate system
$(x,y)$ on the torus such that the metrics has the form $\lambda(x,y)(dx^2 + dy^2)$;
in this coordinates system the form $A$ is %
$(a+ b\cplxi)\frac{\partial}{\partial z}\otimes\frac{\partial}{\partial z}\otimes\frac{\partial}{\partial z}$, where $a, b \in \mathbb{R}$.  Then, one can view the
equations $\calg_2=0$, $\calg_3=0$ on the torus as two PDE on the coefficient
$\lambda $ depending on two parameters $a,b$.
  
One should also mention that in paragraph  \ref{calcul} we invented a calculus
adapted to our problem. The calculations are much easier in this calculus (for
example, in the recent papers \cite{Kr} and \cite{Br-Du-Ea}, which deal with a
priori easier quadratic integrals, most calculations were done using computer
algebra programs; in our paper, everything is done ``by hand'').  We expect
that the same calculus will effectively work in the case of integrals of
higher degree.

{\bf Acknowledgement:} We thank B. Kruglikov for useful discussions, and 
 Deutsche Forschungsgemeinschaft (Priority Program 1154
---  Global Differential Geometry and Research Training Group 1523 --- Quantum and
Gravitational Fields) and FSU Jena for partial financial support.


\

\ifx\undefined\bysame
\newcommand{\bysame}{\leavevmode\hbox to3em{\hrulefill}\,}
\fi

\def\entry#1#2#3#4\par{\bibitem[#1]{#1}
{\textsc{#2 }}{\sl{#3} }#4\par\vskip2pt}

\def\noentry#1#2#3#4\par{}

\def\mathrev#1{{{\bf Math.\ Rev.:\,}{#1}}}

\end{document}